\documentstyle{amsppt}
\magnification=\magstep1
\pagewidth{5.5in}
\vsize=7in 
\widestnumber \key{AAAAAAAAAA}
\topmatter    
\title Some sharp restriction theorems for homogeneous manifolds \endtitle
\author Laura De Carli and Alex Iosevich 
\endauthor
\endtopmatter
\document

\noindent{\bf Introduction}: Let $M$ denote a submanifold of ${\Bbb R}^{n+l}$ of 
codimension $l$. Let $\Cal R$ denote a restriction operator
$${\Cal R}f(\eta) =\int e^{-i\langle x,\,\eta\rangle} f(x)dx,  
\quad  \eta\in M, \quad f\in{\Cal S}({\Bbb R}^{n+l}). 
\tag0.1$$ 

We wish to find an optimal range of exponents $p$ such that
$${||{\Cal R}f||}_{L^2(M,\, d\sigma)} \leq
C_p{||f||}_{L^p({\Bbb R}^{n+l})}, \tag0.2$$
where $d\sigma$ is a compactly supported measure on $M$. 

When $M$ is a codimension one surface in ${\Bbb R}^{n+1}$ with non-vanishing 
Gaussian curvature, the estimate $(0.2)$ is well understood. A celebrated 
result due to Stein and Tomas says that if $M$ has non-vanishing Gaussian  
curvature, then the estimate $(0.2)$ holds with $p=\frac{2(n+2)}{n+4}$. In 
higher dimensions the situation is more complicated. When $M$ is a codimension 
two surface in ${\Bbb R}^{n+2}$ satisfying a non-degeneracy assumption, the
estimate $(0.2)$ holds with $p=\frac{2(n+4)}{n+8}$. (See \cite{MC} and 
Theorem A in the Section 2 below). However, a sharp necessary and sufficient
condition is not currently available.  

It should also be noted that even in codimension one, the more general 
$(L^p, L^q)$ estimates for the restriction operator are not fully understood,
except in the dimension two. (See \cite{St} for a detailed discussion). We 
shall not address this issue here.    

The purpose of this paper is to establish the estimate $(0.2)$ in the case when
$$M=\{(x,x_{k+1},..., x_{k+l}): x_{k+1}={\phi}_1(x),..., x_{k+l}={\phi}_l(x)\},$$
where ${\phi}_j \in C^{\infty}({\Bbb R}^n \backslash 0)$ is homogeneous of degree 
$m_j \ge 1$. 

Let ${\Cal F}[d\sigma]$ denote the Fourier transform of $d\sigma$. By a theorem
of Greenleaf (see \cite{Gr}), the inequality $(0.2)$ holds for 
$\dsize p=\frac{2(l+\gamma)}{2l+\gamma}$ if 
$$\vert{{\Cal F}[d\sigma](R\zeta)}\vert\leq C(1+R)^{-\gamma},
\qquad \zeta \in S^{n+1}. \tag0.3$$  

We shall see below that isotropic Fourier transform estimates do not yield
the sharp restriction theorem in codimension two or higher (see e.g \cite{MC}).
It should be noted that even in codimension one, it is not known whether 
the exponent given by Greenleaf's theorem is sharp. 

The best possible isotropic rate of decay of ${\Cal F}[d\sigma]$ for the  
homogeneous manifold $M$ defined above is $-\frac{n}{\max_j\{m_j\}}$. An 
application of Greenleaf's theorem yields the estimate $(0.2)$ with 
$$p_0=\frac{2(n+lm)}{n+2lm},\tag0.4$$ where $m=\max_j\{m_j\}$.

However, the following homogeneity argument due to Knapp suggests that the
optimal exponent for the estimate $(0.2)$ is    
$$ p_0=\frac{2(m_1+...+m_l+n)}{n+2(m_1+...+m_l)}.\tag0.5$$ 

Indeed, let ${\Cal R}$ denote the restriction operator defined above. Let 
$\hat f_\delta = h$, where $h$ is the characteristic function of a rectangle in 
${\Bbb R}^{n+l}$ with sides of lengths $ (1,\,1,\,\cdots 1,\,C,...,\ C)$, 
$C$ large.

Then  
$$ \Vert{f_\delta}\Vert_p\approx \delta^{(1-1/p)(n+m_1+...+m_l)}\qquad
\text{and} \qquad \Vert{{\Cal R}f_\delta}\Vert_2 \approx\delta^{n/2}.\tag0.6$$

Hence $(0.2)$ can only hold if 
$\dsize p \leq \frac{2(n+m_1+...+m_l)}{n+2(m_1+...+m_l)}$. 

We will establish the estimate $(0.2)$ for a homogeneous manifold $M$, with
the exponent $p_0$ given by Knapp's homogeneity argument, under a variety
of conditions on the level sets of the graphing functions ${\phi}_1$, ${\phi}_2$, 
..., ${\phi}_l$. 
\vskip.125in  

Our main results are the following. The first result gives us a good 
description of $L^2$ restriction theorems for two-dimensional submanifolds
in codimension $2$ given as graphs of homogeneous polynomials. 
(Please see Definition 2.1 and 2.2 below
for the precise description of finite type. Please see Definition 2.8 below
for the description of the order of vanishing along a line.) 

\proclaim{Theorem I} (See Theorem 2.9 in Section 2). 
Let 
$$S=\{(x_1,x_2,x_3,x_4)\ : \ x_3=\phi_1(x_1,x_2), \  x_4=\phi_2(x_1,x_2), \}
\tag0.7$$ where $\phi_1$, $\phi_2\in C ^{\infty}({\Bbb R}^2) $ are homogeneous  
polynomials of degree 
$m_1$ and $m_2$ respectively, $(m_1\ge m_2\ge 2)$. Suppose that there exists a 
non-zero constant $c$ such  that ${{\phi_1(x)}_{|}}_{\{x:\phi_2(x)=1\}} = c $. 
Let $Z_0=\{(x_1,x_2): {\phi}_2(x_1,x_2)=0\}$, $Z_1=\{(x_1,x_2): 
\bigtriangledown {\phi}_2(x_1,x_2)=(0,0)$, and 
$Z_2=\{(x_1,x_2): H{\phi}_2(x_1,x_2)=0\}$, where $H{\phi}_2$ denotes the 
determinant of the hessian matrix of ${\phi}_2$.  
Suppose that the curve $\{x: \phi_2(x) =1\}$ is of finite type 
$m$ at each point of $Z_2\cap \{x: \phi_2(x) =1\}$ and that 
$\phi_2$ vanishes of order $\leq M$ along the lines contained in  
$Z_0 \cup Z_1$.
Then $(0.2)$ holds for every $p<p_0$, where $p_0$ is the sharp exponent 
$$\frac{m_1+m_2+2}{m_1+m_2+1}, \tag0.8$$
if $\dsize m_1+m_2 \ge \max\left\{ 2 M\left(1+\frac{m_1}{m_2}\right),  
2m\right\}$. Furthermore $(0.2)$ holds with 
\roster  
\item
$\dsize
p <\max\left\{ 
\frac{2(1+M(1+\frac{m_1}{m_2})}{1+2M(1+\frac{m_1}{m_2})}, 
\frac{2(m+1)}{2m+1}\right\}$  
if $ \min\{2M(1+\frac{m_1}{m_2}),  2m\}
\leq m_1+m_2 <  \max\{2M(1+\frac{m_1}{m_2}),  2m\}$,
\item
$\dsize
p< \min\left\{ 
\frac{2(1+M(1+\frac{m_1}{m_2})}{1+2M(1+\frac{m_1}{m_2})}, 
\frac{2(m+1)}{2m+1}\right\}$  if 
$4\leq m_1+m_2 < \min\{2M(1+\frac{m_1}{m_2}),  2m\}$.
\endroster
\endproclaim  

The following result answers a question posed by Fulvio Ricci about the
restriction theorems for manifolds given as graphs of quadratic monomials. 
The proof relies on the precise asymptotics of the Fourier transforms of
certain distributions obtained by Shintani (see \cite{TS}, and Lemma 3.3
below). It has been brought to our attention that this result is implied
by a theorem announced in the Bulletin of the AMS by Gerd Mockenhaupt 
(see \cite{GM}). We enclose the proof for reader's convenience, and to 
motivate the related results proved in this paper. (See e.g. Theorem  III 
below).  

\proclaim{Theorem II} (See Theorem 3.1 below). 
Let $S=\{(x,x_{n+1},...,x_{n+l}) \in {\Bbb R}^{n+l}: 
x_{n+1}={\phi}_1(x),..., x_{n+l}={\phi}_l(x)\}$, where $l=\frac{n(n+1)}{2}$, and
the ${\phi}_j$ denote the distinct monomials of degree $2$. Then the estimate
$(0.2)$ holds with the sharp exponent $p_0=\frac{2(n+2)}{2n+3}$. 
\endproclaim 
\vskip.125in 

The following result generalizes Theorem II to manifolds given as joint 
graphs of smooth functions of higher order of homogeneity. The proof relies
on the non-isotropic decay estimates for the associated Fourier transform
of the surface carried measure. The observation that the non-isotropic 
decay estimates are useful to obtain sharp restriction theorems in 
codimension $>1$ is not new. See for example the work of M.Christ (\cite{MC})
and E.Prestini (\cite{P}). 

\proclaim{Theorem III} (See Theorem 3.2 below). 
Let $S$ denote a compact piece of the manifold 
$\{(x,x_{n+1},...,x_{n+l}) \in {\Bbb R}^{n+l}: 
x_{n+1}={\phi}_1(x),..., x_{n+l}={\phi}_l(x)\}$,where $\phi_j \in C^{\infty}
({\Bbb R}^n\slash 0)$
is homogeneous of degree $m \ge 2n$. Suppose that no linear combination 
of the  $\phi_j$'s  vanishes on a subset of postive measure of $S^{n-1}$.
Let ${\Phi}(x)= ({\phi}_1(x),..., {\phi}_l(x))$. 
Suppose that $\Phi(\omega)\not=(0,...,0)$, 
$\omega \in S^{n-1}$. 
Then the estimate $(0.2)$ holds with the 
sharp exponent 
$$p_0=\frac{2(n+lm)}{n+2lm}. \tag0.9$$  
\endproclaim  

Many of our results are based on the non-isotropic decay estimates for the
associated Fourier transform of the surface carried measure. A sample 
result is the following. 
\proclaim{Theorem IV} (See Theorem 1.1B below) 
Let $S=\{(x,x_{n+1},..., x_{n+l}) \in {\Bbb R}^{n+l}: 
x_{n+1}={\phi}_1(x), ..., x_{n+l}={\phi}_l(x)\}$, where each ${\phi}_j \in 
C^{\infty}({\Bbb R}^n\slash 0)$ is homogeneous of degree $m\ge2n$. 
Let ${\Phi}(x)=
({\phi}_1(x),..., {\phi}_l(x))$. Let 
${\Phi}_{\lambda}(x)=\langle \Phi(x), \lambda \rangle$. Let 
$$F(\xi,\lambda)=\int e^{i\langle x,\xi \rangle +{\lambda}_1{\phi}_1(x)+...+
{\lambda}_l {\phi}_l(x))} \psi(x)dx, \tag0.10$$ where $\psi$ is a smooth
cutoff function. Then  
$$ |F(\xi,\lambda)| \leq C\int_{S^{n-1}} 
{|{\Phi}_{\lambda}(\omega)|}^{-\frac{n}{m}} d\omega. \tag0.11$$  
\endproclaim

This paper is organized as follows. In Section 1 we will prove some estimates
related to the decay of the Fourier transform of the surface carried measure. In
particular, we will estimate the decay of the Fourier transform of the surface
carried measure, in any codimension, in terms of the integrability of the  
multiplicative inverses of the graphing functions ${\phi}_j$ restricted to the 
unit sphere. Using this technique we shall also obtain an accurate non-isotropic
estimate for the Fourier transform of the surface carried measure in the case 
when every graphing function has the same homogeneity. 

In Section 2 we will apply the results of Section 1 along with the results of
M. Christ (\cite{MC}), E. Prestini (\cite{P}), 
and a variety of scaling arguments
to obtain a sharp estimate $(0.2)$ with the exponent $p_0$ given by $(0.5)$. 

In Section 3 we will use the non-isotropic estimates from Section 1 to study 
restriction theorems in the case when every graphing function is homogeneous 
of the same degree $m\ge2$. 
\vskip.25in

\head Section 1 \endhead 
\head The decay of the Fourier transform of the surface carried measure \endhead 
\vskip.25in

We will need the following definitions. (See e.g. \cite{Gr}, \cite{MC}).  
\vskip.125in

{\bf Nonvanishing Gaussian curvature:} Let $\Sigma$ be a submanifold of 
${\Bbb R}^{N+1}$ 
of codimension $1$ equipped with a smooth compactly supported measure $d\mu$.
Let $J:\Sigma \to S^{N}$ be the usual Gauss map taking each point on 
$\Sigma$ to the outward  unit normal at that point. We say that $\Sigma$
has everywhere nonvanishing Gaussian curvature if the differential of the 
Gauss map $dJ$ is always nonsingular.
\vskip.125in
\noindent {\bf Strong curvature condition:} Let $S$ be a submanifold of ${\Bbb R}^{n+l}$ 
of codimension $l$ equipped with a smooth compactly supported measure $d\mu$.
Suppose that $S$ is a joint graph of smooth functions $g_1$, $g_2$,...,$g_l$,  
where $g_j:{\Bbb R}^n \rightarrow {\Bbb R}$. 
Let ${\Cal N}_{x_0}(S) $ denote the $l$-dimensional
space of normals to $S$ at a point $x_0$. We say that $S$ satisfies the 
strong curvature condition (SCC) if for all $x_0 \in S$ in some neighborhood of 
the support of $d\mu$, 
$$\det D^2(\nu_1g_1(x)+\nu_2g_2(x)+...+\nu_lg_l(x))\ne 0,\qquad 
\forall \nu\in {\Cal N}_{x_0},\tag1.1$$ 
where $D^2$ denotes the Hessian matrix. 
\vskip.125in

\noindent 
{\bf N-curvature condition:} Let $S\in {\Bbb R}^{n+l}$ be defined as above. We say
that $S$ satisfies the $N$-curvature condition if the rank of the Hessian 
matrix in $(1.1)$ is greater than or equal to $N$ everywhere. 
\vskip.125in

Our main results are the following. 
\vskip.125in

\proclaim{Theorem 1.1A} Let $S=\{(x,x_{n+1},...,x_{n+l})\in {\Bbb R}^{n+l}: 
x_{n+1}={\phi}_1(x),..., x_{n+l}={\phi}_{n+l}(x)\}$, where 
${\phi}_j \in C^{\infty}({\Bbb R}^n \backslash 0)$ is homogeneous of degree
$m_j \ge 2$. Suppose that no ${\phi_i}$ is a constant multiple of any ${\phi_j}$
for $i\ne j$.
Suppose that ${({\phi}_j)}^{-1}\in L^{{\rho}_j}(S^{n-1})$, 
$0<{\rho}_j\leq \frac{n}{m_j}$. Let 
${\mu}$ denote the number of distinct $m_j$'s.  
Suppose that $\frac{n}{m_j} \leq \frac{1}{\mu+1}$ for each $j$. Let 
$\lambda=({\lambda}_1,..., {\lambda}_l)$, and define 
$$ F(\xi, \lambda)=\int e^{i(\langle x,\xi \rangle +{\lambda}_1{\phi}_1(x)+...+
{\lambda}_l{\phi}_l(x))} \psi(x)dx, \tag1.2$$ where $\psi$ is a smooth cutoff
function. 

Then 
$$ |F(\xi,\lambda)| \leq C{(1+|\xi|+|\lambda|)}^{-\rho}, \tag1.3$$ where
$\rho=\min_j \{\rho_j\}$. 
\endproclaim 

\proclaim{Theorem 1.1B} Let $S=\{(x,x_{n+1},..., x_{n+l}) \in {\Bbb R}^{n+l}: 
x_{n+1}={\phi}_1(x), ..., x_{n+l}={\phi}_l(x)\}$, where each ${\phi}_j \in 
C^{\infty}({\Bbb R}^n\slash 0)$ is homogeneous of degree $m\ge2$. 
Let ${\Phi}(x)=
({\phi}_1(x),..., {\phi}_l(x))$. Let 
${\Phi}_{\lambda}(x)=\langle \Phi(x), \lambda \rangle$. Let $F(\xi,\lambda)$ be
defined as in $(1.2)$. Then  
$$ |F(\xi,\lambda)| \leq C\int_{S^{n-1}} 
{|{\Phi}_{\lambda}(\omega)|}^{-\frac{n}{m}} d\omega. \tag1.6$$  
\endproclaim

\proclaim{Lemma 1.2} Let $S\in {\Bbb R}^{n+l}$ and $F(\xi,\lambda)$ be defined 
as above. Let 
$$F_0(\xi,\lambda)=\int e^{i(\langle x,\xi \rangle +{\lambda}_1{\phi}_1(x)+...+
{\lambda}_l{\phi}_l(x))} {\psi}_0(x)dx, \tag1.7$$ where ${\psi}_0$ is a smooth
cutoff function supported away from the origin, and each ${\phi}_j \in 
C^{\infty}({\Bbb R}^n \backslash 0)$ is homogeneous of degree $m\ge2$. Suppose
that 
$$|F_0(\xi,\lambda)| \leq C{(1+|\xi|+|\lambda|)}^{-\gamma}, \ \ \gamma>0. 
\tag1.8$$ 

Then 
$$|F(\xi, \lambda)| \leq C{(1+|\xi|+|\lambda|)}^{-\frac{n}{ \max_j \{m_j\} }}, 
\tag1.9$$ if $\gamma \max_j\{m_j\} \ge n $. 
\endproclaim 

\proclaim{Lemma 1.3} Let $S\in {\Bbb R}^{n+l}$ be defined as above.
Suppose that $S{|}_{\{x_{n+j_1}=1, x_{n+j_2}=1,...,x_{n+j_s}=1\}}$, $1\leq s<l$, 
is an $(n-s)$-dimensional submanifold of codimension $l$, of the hyperplane
$$\{(x,x_{n+1},...,x_{n+l}): x_{n+j_1}=1,..., x_{n+j_s}=1\},$$ satisfying the  
$N$-curvature condition. Let $F_0(\xi,\lambda)$ be defined as in $(1.7)$.  
Then 
$$ |F_{0}(\xi,\lambda)| \leq C{(1+|\xi|+|\lambda|)}^{-\frac{N}{2}}. \tag1.10$$ 
\endproclaim  

\proclaim{Lemma 1.4} Let $S\in {\Bbb R}^{n+l}$ be defined as above 
with $l=2$. Suppose that there exists a constant $c$, such that 
${\phi}_1{|}_{\{x: {\phi}_2(x)=1\}}\equiv c$. Suppose that $m_1\not=m_2$. Suppose
that $H{\phi}_1$, the Hessian matrix of ${\phi}_1$, has rank $\ge 2$ away from 
the origin. Suppose that ${\phi}_1(\omega)\not=0$, $\omega \in S^{n-1}$. 
Let $F_0(\xi,\lambda)$ be defined as in $(1.7)$. Then 
$$ |F_{0}(\xi,\lambda)|\leq C{(1+|\xi|+|\lambda|)}^{-\frac{1}{2}}. \tag1.11$$ 
\endproclaim 

\proclaim{Lemma 1.5} Let $S\in {\Bbb R}^{n+l}$ be defined as in lemma
$1.4$, with $n=l=2$. 
Then the assumptions of Theorem B (see Section 2) 
are satisfied at every $x\ne (0,\,0)$ iff ${\phi}_2(\omega) \not=0$, 
$\omega \in S^1$, and the level set $\{x: {\phi}_2=1 \}$ has non-vanishing 
Gaussian curvature. 
\endproclaim 
\vskip.25in

\demo{Proof of Theorem 1.1A} We shall prove the theorem under the assumption 
that all the $m_j$'s are distinct. The general statement follows by combining  
the terms with the same homogeneity. Recall that 
$$ F(\xi,\lambda)=\int e^{i(\langle x,\xi \rangle +{\lambda}_1{\phi}_1(x)+
...+ {\lambda}_l{\phi}_l(x))} \psi(x)dx, \tag1.12$$ 
where without loss of generality 
$\psi$ is radial, and each ${\phi}_j \in C^{\infty}({\Bbb R}^n)$ is homogeneous\
of degree $m_j \ge 2$. 

Let ${\Omega}=\left\{\omega \in S^{n-1}: 
|{\phi}_1(\omega)|>\frac{1}{|{\lambda}_1|},..., 
|{\phi}_l(\omega)|>\frac{1}{|{\lambda}_l|}\right\}$. 

Let 
$$ F(\xi,\lambda)=\int_{\frac{x}{|x|}\in S^{n-1}\backslash {\Omega}} 
+\int_{ \frac{x}{|x|} \in {\Omega}}=
F_{S^{n-1} \backslash {\Omega}}+F_{{\Omega}}.\tag1.13$$ 

Taking absolute values inside the integral we see that 
$$ |F_{S^{n-1}\backslash {\Omega}}(\xi,\lambda)| \leq 
C|S^{n-1}\backslash {\Omega}|, \tag1.14$$ 
which by Chebyshev's inequality is bounded by
$$ C\max_{j} \left\{{|{\lambda}_j|}^{-\rho_j} 
\int_{S^{n-1}} {|{\phi}_j(\omega)|}^{-\rho_j} d\omega \right\}\leq 
C{|\lambda|}^{-\rho},\tag1.15$$ for $|\lambda|$ large. 

Let ${\psi}_0 \in C^{\infty}_{0}({\Bbb R})$ be supported in the interval 
$[\frac{1}{2}, 4]$, such that ${\psi}_0 \equiv 1$ inside $[1,2]$, and 
$$ \sum_{k=0}^{\infty}{\psi}_0(2^ks)\equiv 1. \tag1.16$$ 

Let 
$$ F^k_{\Omega}(\xi,\lambda)=
\int \int e^{i(r\langle \omega, \xi \rangle+ 
{\lambda}_1r^{m_1}{\phi}_1(\omega)+
...+ {\lambda}_lr^{m_l}{\phi}_l(\omega))} 
{\psi}_0(2^kr)r^{n-1}drd\omega.\tag1.17$$ 

A change of variables sending $r \rightarrow 2^{-k}r$ shows that 
$$F^k_{\Omega}(\xi,\lambda)=2^{-nk}
F^{0}_{\Omega}(2^{-k}\xi, 2^{-m_1k}
{\lambda}_1,..., 2^{-m_lk}{\lambda}_l).\tag1.18$$ 

Let ${\Cal A}=\{1,2,..., l\}$. 

We must estimate 
$$ \sum_{k=0}^{\infty}2^{-nk}|F^{0}_{\Omega}(2^{-nk}\xi, 2^{-m_1k}{\lambda}_1, 
..., 2^{-m_lk}{\lambda}_l)|=
\sum_{{\Cal B} \subset {\Cal A}} I_{{\Cal B}}, \tag1.19$$ where
$$ I_{{\Cal B}}=\sum_{\{{\lambda}_j{\phi}_j(\omega)\leq 2^{m_jk}, \ j\in {\Cal B}, 
\ \ {\lambda}_j{\phi}_j(\omega)>2^{m_jk}, \ j \notin {\Cal B}\}} 
2^{-nk}|F^{0}_{\Omega}(2^{-k}\xi, 2^{-m_1k}{\lambda}_1,..., 2^{-m_lk}{\lambda}_l)|
=\tag1.20$$
$$ \sum_{S_{\{{\Cal B},\leq\}}\bigcap S_{\{{\Cal B},>\}}}  
2^{-nk}|F^{0}_{\Omega}(2^{-k}\xi, 2^{-m_1}k{\lambda}_1,
..., 2^{-m_lk}{\lambda}_l)|.$$ 

To estimate each $I_{{\Cal B}}$ we shall need the fact that the curve 
$(r,r^{m_1},..., r^{m_l})$ has non-vanishing curvature and torsion away from the
origin, so long as all the $m_j$'s are distinct. An elementary van der Corput
type estimate shows that the Fourier transform of the measure carried by this
curve decays of order $-\frac{1}{l+1}$. We shall also use the fact that 
$F^{0}_{\Omega}$ is bounded. More precisely, 

$$ |I_{{\Cal B}}| \leq \int_{\Omega} \left( C\sum_{S_{\{{\Cal B},\leq\}}} 2^{-nk}
\right) +\left(C\sum_{S_{\{{\Cal B},>\}}} 2^{-nk} {\left( \sum_{j \in {\Cal B}} 
2^{-m_jk}|{\phi}_j(\omega){\lambda}_j| \right)}^{-\frac{1}{l+1}} \right)d\omega 
\leq \tag1.21$$
$$ \max_{j \in {\Cal A}} {|{\lambda}_j|}^{-\frac{n}{m_j}} 
\int_{\Omega} {|{\phi}_j(\omega)|}^{-\frac{n}{m_j}}d\omega ,$$ provided that 
$\frac{n}{m_j} \leq \frac{1}{l+1}$, for each $j \in {\Cal A}$. 

Moreover, 
$$ {|{\lambda}_j|}^{-\frac{n}{m_j}} 
\int_{\Omega} {|{\phi}_j(\omega)|}^{-\frac{n}{m_j}}
d\omega={|{\lambda}_j|}^{-\rho_j} \int_{\Omega} {|{\phi}_j(\omega)|}^{-\rho_j}
{|{\phi}_j(\omega){\lambda}_j|}^{\rho_j-\frac{n}{m_j}} d\omega \leq \tag1.22$$
$$ {|{\lambda}_j|}^{-\rho_j} \int_{S^{n-1}} {|{\phi}_j(\omega)|}^{-\rho_j}
d\omega,$$ since $\rho_j \leq \frac{n}{m_j}$. 

Hence, the expression $(1.21)$ is bounded by
$$ \max_{j \in {\Cal A}} {|{\lambda}_j|}^{-\rho_j} \int_{S^{n-1}} 
{|{\phi}_j(\omega)|}^{-\rho_j} d\omega \leq C{(|\lambda|)}^{-\rho}, \tag1.23$$ for
$|\lambda|$ large. 

This completes the proof if $|\xi| \leq C|\lambda|$. However, if this is not
the case, the gradient of the phase function $\langle x, \xi \rangle +
{\lambda}_1{\phi}_1(x)+...+ {\lambda}_l{\phi}_l(x)$ is bounded away from zero, and
an integration by parts argument (see \cite{St}, p.364) shows that 
$|F(\xi,\lambda)| \leq C_N{1+|\xi|+|\lambda|}^{-N}$, for any $N>0$. 
\enddemo
\vskip.125in 

\demo{Proof of Theorem 1.1B} Let $\lambda=t\nu=t({\nu}_1,..., {\nu}_l)$. Let
${\Phi}_{\nu}(x)=\langle \Phi(x), \nu \rangle$. We rewrite $F(\xi,\lambda)$ in
the form
$$ \int e^{i(\langle x,\xi \rangle +t{\Phi}_{\nu}(x))} \psi(x)dx, \tag1.24$$ 
where $\psi$ is a smooth cutoff function. Let  
$\Omega=\{\omega \in S^{n-1}: |{\Phi}_{\nu}(\omega)|>\frac{1}{t}\}$. Let 
$$ F(\xi,\lambda)=\int_{\frac{x}{|x|}\in S^{n-1} \backslash \Omega} +
\int_{\frac{x}{|x|} \in \Omega}=F_{S^{n-1}\backslash \Omega}+ 
F_{\Omega}. \tag1.25$$ 

Taking the absolute values inside the integral, we see that
$$ |F_{S^{n-1}\backslash \Omega} (\xi,\lambda)| \leq 
C|S^{n-1}\backslash \Omega|=C|\{\omega \in S^{n-1}: 
|{\Phi}_{\nu}^{-1}(\omega)|>t\}|, \tag1.26$$ which by Chebyshev's 
inequality is bounded by 
$$ t^{-\frac{n}{m}} \int {|{\Phi}_{\nu}(\omega)|}^{-\frac{n}{m}} d\omega. 
\tag1.27$$  

Let ${\psi}_0\in C^{\infty}_0({\Bbb R})$ be supported in the interval 
$[\frac{1}{2}, 4]$, such that ${\psi}_0 \equiv 1$ inside $[1,2]$, and 
$$ \sum_{k=0}^{\infty} {\psi}_0(2^ks) \equiv 1. \tag1.28$$ 

Let 
$$F^{k}_{\Omega}(\xi,\lambda)=\int \int e^{i(r\langle \omega, \xi \rangle +
tr^m{\Phi}_{\nu}(\omega))} r^{n-1} {\psi}_0(2^kr) d\omega dr. \tag1.29$$ 

A change of variables sending $r \rightarrow 2^{-k}r$ shows that 
$F^{k}_{\Omega}(\xi,\lambda)=2^{-nk}F^{0}_{\Omega}(2^{-k}\xi, 2^{-mk}\lambda).$
We must estimate 
$$\sum_{k=0}^{\infty}2^{-nk}|F^{0}_{\Omega}(2^{-k}\xi, 2^{-mk}\lambda)|=
\sum_{\{t{\Phi}_{\nu}(\omega) \leq 2^{mk}\}} 
+\sum_{\{t{\Phi}_{\nu}(\omega)>2^{mk}\}}=I+II. \tag1.30$$ 

To estimate $II$ we shall use the fact that away from zero the Fourier transform
of the measure supported on the curve $(r,r^m)$ decays of order $-\frac{1}{2}$. To
estimate $I$ we shall just use the fact that $F^{0}_{\Omega}$ is bounded. More 
precisely, 
$$ |I| \leq C \int_{\Omega}  
\sum_{\{t{\Phi}_{\nu}(\omega)\leq 2^{mk}\}} 2^{-nk} d\omega \leq
C t^{-\frac{n}{m}} \int {|{\Phi}_{\nu}(\omega)|}^{-\frac{n}{m}} d\omega, 
\tag1.31$$  and
$$ |II| \leq C \int_{\Omega}  \sum_{\{t{\Phi}_{\nu}(\omega)>2^{mk}\}} 
2^{-nk} {|2^{-mk}{\Phi}_{\nu}(\omega)t|}^{-\frac{1}{2}}  d\omega \leq 
Ct^{-\frac{n}{m}} \int_{\Omega} 
{|{\Phi}_{\nu}(\omega)|}^{-\frac{n}{m}} d\omega, \tag1.32$$
as long as $m\ge2n$. 

This completes the proof if $|\xi| \leq C|\lambda|$. However, if this is not the
case, the gradient of the phase function 
$\langle x, \xi \rangle +{\lambda}_1{\phi}_1(x)+...+ {\lambda}_l{\phi}_l(x)$ is
bounded away from the origin, and an integration by parts argument  
(see \cite{St}, p.364) shows that 
$|F(\xi,\lambda)| \leq C_N{(1+|\xi|+|\lambda|)}^{-N}$, for any $N>0$. 
\enddemo 

\demo{Proof of Lemma 1.2} Let 
$$F_k(\xi,\lambda)=\int e^{i(\langle x,\xi \rangle +{\lambda}_1{\phi}_1(x)+ 
...+ {\lambda}_l{\phi}_l(x))} {\psi}_0(2^kx)dx. \tag1.33$$ 

A change of variables sending $x \rightarrow 2^{-k}x$ shows that 
$$F_k(\xi,\lambda)=2^{-nk}F_0(2^{-k}\xi, 2^{-m_1k}{\lambda}_1,..., 
2^{-m_lk}{\lambda}_l). \tag1.34$$ 

Let ${\Cal A}=\{1,2,..., l\}$. We must estimate 
$$ \sum_{k=0}^{\infty} 2^{-nk} |F_0(2^{-k}\xi, 2^{-m_1k}{\lambda}_1,..., 
2^{-m_lk}{\lambda}_l)|=\sum_{{\Cal B}\subset {\Cal A}} I_{{\Cal B}}, \tag1.35$$ 
where
$$ I_{{\Cal B}}=\sum_{\{|{\lambda}_j|\leq 2^{m_jk}, \ j\in {\Cal B}, \ \ 
|{\lambda}_j|>2^{m_jk}, \ j\notin {\Cal B}\}} 2^{-nk}
|F_0(2^{-k}\xi, 2^{-m_1k}{\lambda}_1,..., 2^{-m_lk}{\lambda}_l)|= \tag1.36$$
$$\sum_{S_{\{{\Cal B},\leq\}} \bigcap S_{\{{\Cal B}, >\}}}2^{-nk} 
|F_0(2^{-k}\xi, 2^{-m_1k}{\lambda}_1,..., 2^{-m_lk}{\lambda}_l)|.$$ 

Using the assumed decay of $F_0$, and the fact that, in particular, $F_0$ is
bounded, we get 
$$ |I_{{\Cal B}}| \leq C_1\sum_{S_{\{{\Cal B},\leq \}}} 2^{-nk} + 
C_2 \sum_{S_{\{{\Cal B},>\}}} 2^{-nk} {\left( \sum_{j \notin {\Cal B}} 
2^{-m_jk}|{\lambda}_j| \right)}^{-\gamma} \leq \tag1.37$$
$$ C\max_{j \in {\Cal A}} {|{\lambda}_j|}^{-\frac{n}{m_j}} \leq 
C {|\lambda|}^{-\frac{n}{\max_j\{m_j\}}}, \tag1.38$$ as long as 
$\max_j\{m_j\} \leq n\gamma$.

This completes the proof if $|\xi|\leq C|\lambda|$. However, if this is not the
case, the gradient of the phase function $\langle x,\xi \rangle +
{\lambda}_1{\phi}_1(x)+...+ {\lambda}_l{\phi}_l(x)$ is bounded away from the 
origin. An integration by parts argument, (see \cite{St}, p.364), shows that
$|F(\xi,\lambda)| \leq C_N{(1+|\xi|+|\lambda|)}^{-N}$, for any $N>0$. 
\enddemo

\demo{Proof of lemma 1.3} Assume that $x_{n+j_1} = x_{n+1}$, 
$\cdots x_{n+j_s} = x_{n+s}$.
In what follows we will denote  $x\in {\Bbb R}^n$ by 
$(x^{\prime},\,x^{\prime \prime})$, 
with $x^{\prime} \in {\Bbb R}^{s}$, $x^{\prime \prime} \in {\Bbb R}^{n-s}$ and  
by $J_{x'}(f)$ (resp. $J_{x''}(f)$) and $D_{x'}^2g$ 
(resp. $D^2_{x''}g$) the Jacobian 
of a function $f(x',\,x'')$ and the Hessian matrix of a 
function $g(x',\,x'')$ computed with respect to the $x'$ (resp. $x''$) variables. 
Let $\phi=(\phi_1,..., \phi_s)$. 

Take  $P$ on the support of $d\sigma$ such that $\phi(P)=(1,1,...,1).$ 
Since $\dsize S_{\vert_{\{x_{n+1}=1, \cdots x_{n+s}=1\}}}$ is by assumption 
a submanifold of codimension $l$ of the hyperplane 
$\{(x, x_{n+1},\cdots, x_{n+l}) : x_{n+1}=1, \cdots x_{n+s} = 1 \}$,  
the Jacobian of the function 
$\phi$, $J(\phi)$, has rank $s$ at $P$.

There is no loss of generality if we assume  that $J_{x'}\phi(P)$ is the 
identity in the space of $s\times s$ matrices, and that $J_{x''}\phi(P)=0$ 

By the implicit function theorem there exists a smooth function 
$\psi(x''):{\Bbb R}^{n-s}\to {\Bbb R}^{s}$ such that 
$\phi(\psi(x''), x'') \equiv 1$   in a 
neighborhood of $P''$.

An application of the chain rule yields:
\roster
\item $J\phi(\psi(x^{\prime \prime}), x^{\prime \prime})=
J_{x'}\phi(\psi(x''),x'')\times J\psi(x'') +J_{x''}\phi(\psi(x''),x'')\equiv 0 $,
which implies that $J\psi(P'') =0$, 
\vskip.1 in
\item
$
{\frac{\partial}{\partial}}_{x_{s+j}} J\psi(P'') + 
{\frac{\partial}{\partial}}_{x_{s+j}}J_{x''}\phi(P) =0
$
for every $j\leq n+s$, and 
\vskip.1 in
\item
${{D^2(\phi_{s+j})(\psi(x''),\,x'')}_\vert}_{x''=P''} = \
\sum_{k=1}^{s}{\frac{\partial}{\partial}}_{x_k}
{\phi_{s+j}(\psi(x''),\,x'')_\vert}_{x''=P''}D^2(\psi_k)(P'')$ 

$ + D^2_{x''}\phi_{s+j}(\psi(x''),\,x'')$
for every $j\leq n-s$.
\endroster
>From (2) we have that, for every  $k\leq s$, the Hessian matrix of $\psi_k$ at 
$P''$, $D^2\psi_k(P'')$, is $-D^2_{x''}\phi_k(P'')$, 
and from (3) that  $D^2_{x''}\phi_{s+j}(P)$ 
can be written as a linear combination of the Hessian 
matrices of the functions $\psi_k$  and of the Hessian matrix of  
$\phi_{s+j}(\psi(x''),\,x'')$ at $P''$.

Let $\Phi =\langle x,\, \xi\rangle +\lambda_1\phi_1+\cdots \lambda_l\phi_l $
denote the phase function of  $F$   as in (1.2). 
By the above remark, $D^2_{x''}\Phi(P)$  can be written 
as a linear combination of the Hessian matrices of the functions  
$\psi_k(x'')$ and  the function $\phi_{s+j}(\psi(x''),\,x'')$ at $P''$. Since 
$\dsize S_{\vert_{\{x_{n+1}=1, \cdots x_{n+s} = 1 \}}}$ 
satisfies the $N$-curvature condition, 
the rank of every linear combination of the above 
matrices  is $N$  for every $\lambda\ne 0$. This shows that the rank of the 
Hessian matrix of $\Phi$ is $\ge N$ and hence that (1.9) holds.
\enddemo 

\demo{Proof of Lemma 1.4} 
A theorem of Littman (\cite{L}) says that if a surface in codimension one has 
at least $k$ non-vanishing principal curvatures, then the Fourier transform of
the surface carried measure decays of order $-\frac{k}{2}$. The proof of that 
theorem shows that $F_0$ has the required decay if the rank of the Hessian 
matrix of the phase function
$\Phi(x)=\langle x,\xi\rangle+ \lambda_1 \phi_1(x)+\lambda_2\phi_2(x)$
is $\ge 1$ for every $x$ on the support of $\chi$ and for every 
$(\lambda_1,\,\lambda_2) \ne (0,\,0)$. 

We  observe now that  $\phi_1 (x)\equiv c\phi_2^{\frac{m_1}{m_2}}(x)$. Indeed, 
if $\phi_2(x_0) \ne 0$, then 
$\dsize \phi_2\left(\frac{x_0}{\phi_2(x_0)^{\frac{1}{m_2}}}\right)= 1$, and 
since $\phi_1{|}_{\{\phi_2(x)=1\}} \equiv c$, then
$\dsize \phi_1\left(\frac{x_0}{\phi_2(x_0)^{\frac{1}{m_2}}}\right)= 
\frac{\phi_1(x_0)}{\phi_2(x_0)^{\frac{m_1}{m_2}}} =c $.
Then $\phi_1(x_0) =c\phi_2(x_0)^{\frac{m_1}{m_2}}$ 
for every $x_0\in {\Bbb R}^n$  such that
$\phi_2(x_0)\ne 0$. But  the argument that we have just used shows
that    $\phi_2(x_0)\ne  0$ iff $\phi_1(x_0)\ne  0$. Hence 
$\phi_1(x) =c\phi_2(x)^{\frac{m_1}{m_2}}$ for every $x \in {\Bbb R}^n$.
\bigskip
By the chain rule,
$$D_{i,j}\phi_1(x)=c\frac{m_1}{m_2} 
D_{i,j}\phi_2(x)\phi_2(x)^{\frac{m_1}{m_2}-1}+ c\frac{m_1}{m_2}
\left(\frac{m_1}{m_2}-1\right)D_i(\phi_2)
D_j(\phi_2)\phi_2(x)^{\frac{m_1}{m_2}-2}, \tag1.39$$
for $i,j, \leq n$. Then  
$$D^2(\Phi)(x)= 
\left(\lambda_2+ c\lambda_1\frac{m_1}{m_2}
\phi_2(x)^{\frac{m_1}{m_2}-1}\right) D^2(\phi_2)(x)+ \tag1.40$$
$$c\lambda_1\frac{m_1}{m_2}\left(\frac{m_1}{m_2}-1\right)
\phi_2(x)^{\frac{m_1}{m_2}-2}
\left(D_{1}(\phi_2)(x)\bigtriangledown(\phi_2)(x)\cdots 
D_{n}(\phi_2)(x)\bigtriangledown(\phi_2)(x)\right).$$ 

Let $P$ be a point on the support of $\chi$. 
If we set $\lambda_2+c\lambda_1\frac{m_1}{m_2}\phi_2(P)^{\frac{m_1}{m_2}-1}=0$, 
we observe that $D^2(\Phi)(P)$ equals a matrix 
whose rank is $1$ if and only if 
$\bigtriangledown(\phi_2)(P)\ne ( 0,\cdots, 0)$,  
$\phi_2(P)^{\frac{m_1}{m_2}-2}\ne 0$, and zero elsewhere. Since we assumed 
that $\phi_2$ doesn't vanish away from the origin,  
Euler's homogeneity relations guarantee that 
$\bigtriangledown(\phi_2)(P)$ doesn't vanish. Consequently, the rank of   
$D^2(\Phi)(P)$ is at most $1$.

To show that $D^2(\Phi)(P)$ cannot be zero,  
we observe that if this were the case we would have
\vskip .1 in
$$-\frac{\left(\lambda_2+c\lambda_1\frac{m_1}{m_2}
\phi_2(P)^{\frac{m_1}{m_2}-1} \right)}
{c\lambda_1\frac{m_1}{m_2}\left(\frac{m_1}{m_2}-1\right)\phi_2(P)^
{\frac{m_1}{m_2}-2}}
D^2\phi_2(P) =
\left(D_{1}(\phi_2)(P)\bigtriangledown(\phi_2)(P)... 
D_{n}(\phi_2)(P)\bigtriangledown(\phi_2)(P)\right). \tag1.41$$

The coefficient which multiplies $D^2\phi_2(P)$ 
cannot be zero because the matrix on the 
right-hand side has rank one. On the other hand, the  
matrix on the left-hand side has rank $\ge 2$ by assumption, 
hence the equality in $(1.41)$ can never hold. 

This concludes the proof of the theorem.
\enddemo

\demo{Proof of Lemma 1.5}
After perhaps rotating and dilating the coordinates, we can   
work in a neighborhood of the point  $(0,1)$. 

Since there exists a constant $c$ such that 
$\phi_1 {|}_{\{x: {\phi}_2=1\}} \equiv c$, 
$\phi_1\equiv c\phi_2^{\frac{m_2}{m_2}}$. (See proof of Lemma 1.4).   

By the chain rule
$D_{i,j}\phi_1(x_1,x_2)=$
$$\frac{m_1}{m_2} 
D_{i,j}\phi_2(x_1,x_2)\phi_2(x_1,x_2)^{\frac{m_1}{m_2}-1}+
\frac{m_1}{m_2}\left(
\frac{m_1}{m_2}-1\right)
D_i(\phi_2)D_j(\phi_2)\phi_2(x_1,x_2)^{\frac{m_1}{m_2}-2}, \tag1.42$$
for $\ i,j, \leq 2$.
\vskip.125in 

Observe that by Euler's homogeneity relations
\roster{}
\item
$D_1(\phi_2)(0,1)=m_2\phi_2(0,1)$,
\item
$D_2(\phi_2)(0,1)=\frac{1}{m_2-1}D_{12}\phi_2(0,1)$, and
\item
$D_{11}(\phi_2)(0,1) = m_2(m_2-1)\phi_2(0,1)$.
\endroster
\vskip.125in

In order to  show that the sufficient condition of Theorem 1.1B is 
verified, we must show that the  determinant of 
the matrix  $J(x,y),$ whose rows are 
$\bigtriangledown(\frac{x_1^2}{2}D_{11}\phi_1(0,1) 
$\par\noindent
$+ \frac{x_2^2}{2}D_{22}\phi_1(0,1) + x_1x_2D_{12}\phi_1(0,1) )$ and  
$\bigtriangledown(\frac{x_1^2}{2}D_{11}\phi_2(0,1) +
\frac{x_2^2}{2}D_{22}\phi_2(0,1)+x_1x_2D_{12}\phi_2(0,1)),$  
is not a square.  
A direct computations shows that the discriminant 
of the determinant of $J(x_1,x_2)$ is
$$\frac{
m_1^2 (\frac{m_1}{m_2}-1)^2}
{(m_2-1)^2} (\phi_2(0,1))^{2(\frac{m_1}{m_2}-1)}[\det(H(\phi_2))(0,1)]^2, \tag1.43$$
which doesn't vanish by the assumptions on $\phi_2$. Hence, $J^{-a}$ is integrable
for $a<1$. This concludes the proof of the lemma. 
\enddemo 
\vskip.75in
\head Section 2 \endhead
\head Restriction theorems-Scaling \endhead 
\vskip.125in

We will need the following results. 
\vskip.125in
\proclaim{Theorem A (M. Christ (\cite{MC}))} Let 
$$S=\{(x,x_{n+1},x_{n+2})\in {\Bbb R}^{n+2}: x_{n+1}=g_1(x), \ \  
x_{n+2}=g_2(x)\}, \tag2.1$$ where $g_j\in C^{\infty}({\Bbb R}^n)$. 

Suppose that $(0.2)$ holds with $p_0=\frac{2(n+4)}{n+8}$. Then for every 
$k \leq \frac{n}{2}$ and every $\theta$, the expression  
$$ {\left(\frac{d}{d\theta}\right)}^k 
\left(det D^2(cos(\theta)g_1(x)+sin(\theta)g_2(x))\right)
\tag2.2$$ does not vanish.  
\endproclaim

\proclaim{Theorem B (E. Prestini (\cite{P}))} Let 
$$ S=\{(x,x_{n+1},..., x_{2n}) \in {\Bbb R}^{n+n}: x_{n+1}=g_1(x), ..., 
x_{2n}=g_n(x)\}, \tag2.3$$ where $g_j \in C^{\infty}({\Bbb R}^n)$. Let 
$G_j$ denote the quadratic part of the Taylor expansion of $g_j$. 
Suppose that the vectors 
$\left\{\bigtriangledown\left( \frac{\partial}{\partial x_i}
G_j\right)(0)\right\}_{{i\leq n}\atop{j\leq n}}$ span 
${\Bbb R}^n$. Let 
$J(x)$ denote the determinant of the matrix 
$({\bigtriangledown}_x G_1,..., {\bigtriangledown}_x G_n)$. Suppose that 
$J^{-a} \in L^1(S^{n-1})$, for any $a<1$. Then $(0.2)$ holds with
$p_0=\frac{6}{5}$. 
\endproclaim 

Note that when $n=2$, the assumptions of Theorem B are equivalent to the 
necessary condition $(2.2)$ in Theorem A. In particular, the conditions of 
Theorem B are necessary and sufficient in that case. 

Before stating our main results, we need to introduce the following definitions.

\vskip.125in

\proclaim{Definition 2.1} Let $\Gamma:I\rightarrow {\Bbb R}^2$, where $I$ is a 
compact interval in ${\Bbb R}$, and $\Gamma$ is smooth. We say that $\Gamma$ is 
finite type if $\langle (\Gamma(x)-\Gamma(x_0)), \mu \rangle$ 
does not vanish of infinite
order for any $x_0\in I$, and any unit vector $\mu$. 
\endproclaim 

We will also need a more precise definition to specify the 
order of vanishing at each point. Let $a_0$ denote a point
in the compact interval $I$. We can always find a smooth function $\gamma$,
such that in a small neighborhood of $a_0$, $\Gamma(s)=(s,\gamma(s))$, where
$s\in I$. 

\proclaim{Definition 2.2} Let $\Gamma$ be defined as before. Let 
$\Gamma(s)=(s,\gamma(s))$ in a small neighborhood of $a_0$. We say that 
$\Gamma$ is finite type $m$ at $a_0$ if ${\gamma}^{(k)}(a_0)=0$ for 
$1\leq k<m$, and ${\gamma}^{(m)}(a_0)\not=0$. \endproclaim  

Our main results are the following.

\proclaim{Lemma 2.3}   
Let $S=\{(x,\,x_{n+1},\cdots,\,x_{n+l}\in {\Bbb R}^{n+l}: 
x_{n+1}=\phi_1(x), \cdots,\, x_{n+l}=\phi_n(x)\}$,  
where each $\phi_j\in{\Cal C}^\infty ({\Bbb R}^n\slash\{0\})$, 
homogeneous of degree $m_j\ge 1$.
Let $d\sigma$ denote  a compactly supported smooth measure on $S$, and let 
$d\sigma_0 =\chi(x)d\sigma$,  where $\chi(x)$ 
is a smooth cutoff  function supported away from 
the origin. Let $Tf(x)=f*\widehat{ d\sigma}$ and let 
$T_0f(x)=f*\widehat d\sigma_0$.
Suppose that $T_0:L^{q_0}({\Bbb R}^{n+l})\to 
L^{q'_0}({\Bbb R}^{n+l})$ is a bounded operator. Then 
$T:L^{p_0(m_1,\cdots m_l)}({\Bbb R}^{n+l}) \to  
L^{p'_0(m_1,\cdots m_l)}({\Bbb R}^{n+l})$ is a bounded 
operator, where 
$$p_0(m_1,\cdots,m_l) = 
\frac{2(n+m_1+\cdots +m_l)}{n+2(m_1+\cdots +m_l)}, \tag2.4$$
as long as $p_0(m_1,\cdots ,m_l)\leq q_0$.
\endproclaim 

\proclaim{Theorem 2.4} Let 
$$S=\{(x,x_{n+1},..., x_{n+l}): x_{n+1}=
{\phi}_1(x), ..., x_{n+l}={\phi}_l(x)\}, \tag2.5$$
where ${\phi}_j \in C^{\infty}({\Bbb R}^n \backslash 0)$ is homogeneous of degree
$m_j\ge 2$. Suppose that $S$ satisfies the assumptions of Theorem $1.1$A.  
Then the estimate $(0.2)$ holds with  
$$ p_0=\frac{2(l+\rho)}{2l+\rho}. \tag2.6$$   
\endproclaim 

\proclaim{Theorem 2.5} Let 
$$S=\{(x,x_{n+1},..., x_{n+l}): x_{n+1}=
{\phi}_1(x), ..., x_{n+l}={\phi}_l(x)\}, \tag2.7$$
where ${\phi}_j \in C^{\infty}({\Bbb R}^n \backslash 0)$ is homogeneous of degree
$m_j \ge 2$. Suppose that $S$ satisfies the assumptions of the Lemma 1.3. Then 
the estimate $(0.2)$ holds with the sharp exponent 
$$ p_0=\frac{2(n+m_1+...+ m_l)}{n+2(m_1+...+m_l)}, \tag2.8$$ provided that 
$(m_1+... +m_l) \ge \frac{2ln}{N}$. 
\endproclaim

\proclaim{Theorem 2.6} Let 
$$ S=\{(x,x_{n+1}, x_{n+2}): 
x_{n+1}={\phi}_1(x), x_{n+2}={\phi}_2(x)\}, \ n>2, \tag2.9$$ where 
${\phi}_j \in C^{\infty}({\Bbb R}^n \backslash 0)$ is homogeneous of degree 
$m_j \ge 2$. Suppose that there exists a non-zero constant $c$, such that 
$$ {\phi}_1 {|}_{\{x: {\phi}_2(x)=1\}} \equiv c.$$ 

Suppose that ${\phi}_1(\omega)\not=0$, $\omega \in S^{n-1}$.  
If the rank of the Hessian matrix of ${\phi}_1$ is $\ge 2$, then 
$(0.2)$ holds with $p_0=\frac{2(m_1+m_2+n)}{n+2(m_1+m_2)}$, provided that 
$m_1+m_2 \ge 4n$. 
\endproclaim  

In order to introduce the Theorem 2.9 below, we need the following result which 
was stated and proved by the second author in (\cite{I2}). 

\proclaim{Lemma 2.7} Let $Z_0=\{(x_1,x_2): P(x_1,x_2)=0\}$. Let 
$Z_1=\{(x_1,x_2): \bigtriangledown P(x_1,x_2)=(0,0)\}$. Let 
$Z_2=\{(x_1,x_2): HP(x_1,x_2)=0\},$ where $HP(x_1,x_2)$ denotes the determinant 
of the Hessian matrix of $P$. Then for each $j=0,1,2$,  
$Z_j=\{(0,0)\} \bigcup {\bigcup}^{N_j}_{k=1}L_k$, where each $L_k$ is a line
through the origin, and $N_j<\infty$. Moreover, $Z_1=Z_0\bigcap Z_2$. 
\endproclaim 

\demo{Proof of Lemma 2.7} Let $P_j$ denote the partial derivative of $P$ with
respect to $x_j$. Since $P$ is homogeneous of degree $m$, $P_j$ is homogeneous
of degree $m-1$, and $HP(x_1,x_2)$ is homogeneous of degree $2(m-2)$. By   
homogeneity, if $Z_j$ contains a point $(x_1,x_2)$, it also contains a line 
through the origin containing that point. Since $P$ is a polynomial, there can
be at most a finite number of such lines. This proves the first assertion of
the lemma. 

By the Euler homogeneity relations, 
$$mP(x_1,x_2)=x_1P_1(x_1,x_2)+x_2P_2(x_1,x_2),$$  
$$(m-1)P_1(x_1,x_2)=x_1P_{11}(x_1,x_2)+x_2P_{12}(x_1,x_2), \tag2.10$$ and 
$$(m-1)P_2(x_1,x_2)=x_1P_{21}(x_1,x_2)+x_2P_{22}(x_1,x_2),$$ where 
the $\{P_{jk}\}$ denote the second partial derivatives. Hence, $Z_0 \subset
Z_1$. If we write the equations for $P_1$ and $P_2$ in matrix form we see
that $(m-1)\bigtriangleup P(x_1,x_2)$ is obtained by applying the Hessian 
matrix of $P$ to the vector $(x_1,x_2)$. Hence, $Z_2\subset Z_1$. Putting
these observations together we see that $Z_0\bigcap Z_2\subset Z_1$. 

Suppose that both $P$ and $HP$ vanish along a line through the origin, which
without loss of generality we take to be the $x_1$-axis. Then
$mP(x_1,0)=x_1P_1(x_1,0)$. This implies that $P_1(x_1,0)=0$. 
Also, $(m-1)P_1(x_1,0)=x_1P_{11}(x_1,0)$. This implies that $P_{11}(x_1,0)=0$.
$P_2(x_1,0)=x_1P_{12}(x_1,0)$. By assumption, 
$$HP(x_1,0)=
P_{11}(x_1,0)P_{22}(x_1,0)-P^2_{12}(x_1,0)=0. \tag2.11$$ 
Since $P_{11}(x_1,0)=0$, we must
conclude that $P_{12}(x_1,0)=0,$ which implies that $P_2(x_1,0)=0$. This proves
that $\bigtriangledown P(x_1,0)=(0,0)$ and hence that $Z_1\subset 
Z_0\bigcap Z_2$. This completes the proof of the lemma. 
\enddemo 
\bigskip
We shall need the following definition:
\vskip.125 in
\noindent
{\bf Definition 2.8} Let $f\in C ^{\infty}({\Bbb R}^2)$. We say
that $f$ vanishes of order $M$ along the line $L=\{(x_1,\, x_2) \ : \ 
x_1=s_1 t, \, x_2 = s_2 t,\ t \in {\Bbb R}\}$ 
if $M$ is the largest positive integer so that 
$f(s_1 t,\,s_2 t)= t^M g(t)$, where $g\in C ^{\infty}({\Bbb R}^2)$ 
is allowed to vanish only at the origin.
\vskip.25in 

\proclaim{Theorem 2.9} 
Let $S=\{(x_1,x_2,x_3,x_4)\ : \ x_3=\phi_1(x_1,x_2), \  x_4=\phi_2(x_1,x_2), \}$
where $\phi_1$, $\phi_2\in C ^{\infty}({\Bbb R}^2) $ are homogeneous  
polynomials of degree 
$m_1$ and $m_2$ respectively, $(m_1\ge m_2\ge 2)$. Suppose that there exists a 
non-zero constant $c$ such  that ${{\phi_1(x)}_{|}}_{\{x:\phi_2(x)=1\}} = c $. 
Let $Z_0$, $Z_1$, $Z_2$   be defined as in Lemma 2.7 
with respect to $\phi_2$. Suppose that  the curve $\{x: \phi_2(x) =1\}$ is of 
finite type
$m$ at each point of $Z_2\cap \{x: \phi_2(x) =1\}$ and that 
$\phi_2$ vanishes of order $\leq M$ along the lines contained in  
$Z_0 \cup Z_1$.
Then $(0.2)$ holds for every $p<p_0$, where $p_0$ is the sharp exponent 
$$\frac{m_1+m_2+2}{m_1+m_2+1}, \tag2.12$$
if $\dsize m_1+m_2 \ge \max\left\{ 2 M\left(1+\frac{m_1}{m_2}\right),  
2m\right\}$. Furthermore $(0.2)$ holds with 
\roster  
\item
$\dsize
p <\max\left\{ 
\frac{2(1+M(1+\frac{m_1}{m_2})}{1+2M(1+\frac{m_1}{m_2})}, 
\frac{2(m+1)}{2m+1}\right\}$  
if $ \min\{2M(1+\frac{m_1}{m_2}),  2m\}
\leq m_1+m_2 <  \max\{2M(1+\frac{m_1}{m_2}),  2m\}$,
\item
$\dsize
p< \min\left\{ 
\frac{2(1+M(1+\frac{m_1}{m_2})}{1+2M(1+\frac{m_1}{m_2})}, 
\frac{2(m+1)}{2m+1}\right\}$  if 
$4\leq m_1+m_2 < \min\{2M(1+\frac{m_1}{m_2}),  2m\}$.
\endroster
\endproclaim

\demo{Proof of Lemma 2.3} Let 
$$ \widehat {d\sigma}(\xi,\,\lambda) =
\int_{  {\Bbb R}^n}
e^{i(
\langle x,\,\xi\rangle +\lambda_1\phi_1(x)+\cdots +\lambda_l\phi_l(x)
)}
\psi(x)dx, \tag2.13$$
where $\psi(x) $ is a cutoff function. 
Let $\rho$ be a cutoff function supported in the 
interval $(1,\,2)$ such that $\sum_{j=0}^{+\infty} \rho(2^jt) =1$ 
for every $t$, and let
$$\widehat{d\sigma_j}( \xi, \lambda) =
\int_{  {\Bbb R}^n}
e^{i(
\langle x,\,\xi\rangle +\lambda_1\phi_1(x)+\cdots +\lambda_l\phi_l(x)
)}
\rho(2^jt)dx. \tag2.14$$

If we make the change of variables sending $x \to 2^{-j}x$ we can write:
$$ \widehat{d\sigma_j}(\xi,\,\lambda)= 
2^{-jn}  \widehat{d\sigma_0}(2^{-j}\xi, \,2^{-m_1j}\lambda_1,\cdots  
2^{-m_lj}\lambda_l). \tag2.15$$

Let $\tau_j$ denote the nonisotropic dilation
$$ \tau_j f(\xi,\,\lambda_1,\cdots,\,\lambda_l)=
f(2^j\xi,\,2^{m_1j}\lambda_1,\,\cdots ,\, 2^{m_lj}\lambda_l). \tag2.16$$

Then
$$ \widehat{d\sigma_j}*f= 
2^{-nj}( \tau_{-j}  \widehat{d\sigma_0}*f)(\xi,\,\lambda)=
2^{-nj}( \tau_{-j}  \widehat{d\sigma_0}*( \tau_{-j} 
\tau_jf))(\xi,\,\lambda). \tag2.17$$

A change of variables shows that 
$ \tau_{-j}  \widehat{d\sigma_0}*(\tau_{-j} \tau_jf)=
2^{j(n+m_1+\cdots m_l)} \tau_{-j}(\widehat{d\sigma_0}*\tau_jf).$

It follows that 
$${||{  \widehat{d\sigma_j}*f}||}_{q'_0}= 
2^{-nj+ j(n+m_1+\cdots m_l)}
{||{ \tau_{-j}(\widehat{d\sigma_0}*\tau_jf)}||}_{q'_0}= \tag2.18$$
$$=2^{j(m_1+\cdots m_l)}
2^{\frac{j}{q'_0} (n+m_1+\cdots m_l)}
{||{  \widehat{d\sigma_0}*( \tau_jf)}||}_{q'_0}\leq$$
$$\leq C 2^{j(m_1+\cdots m_l)}
2^{\frac{j}{q'_0} (n+m_1+\cdots m_l)}
{||{ \tau_jf}||}_{q_0}= 
C2^{j(m_1+\cdots m_l)}
2^{j\left(\frac{1}{q'_0}-\frac{1}{q_0}\right) (n+m_1+\cdots m_l)}
{||{f}||}_{q_0}.$$

The series $\dsize
\sum_{j=1}^{\infty}{||{\widehat{d\sigma_j}*f}||}_{q'_0}$ converges,  provided that 
$$m_1+\cdots m_l +
\left(\frac{1}{q'_0}-\frac{1}{q_0}\right) (n+m_1+\cdots m_l)<0, \tag2.19$$
which  yields
$$q_0< \frac{ 2(n+m_1+\cdots m_l)}{n+2(m_1+\cdots m_l)}. \tag2.20$$

This concludes the proof of the lemma.
\enddemo

\demo{Proof of Theorem 2.4} The application of Greenleaf's theorem (see (0.3)
above) yields the result. 
\enddemo

\demo{Proof of Theorem 2.5} Let $F_0$ be defined as in $(1.7)$. By Lemma 1.3 
$$ |F_0(\xi,\lambda)| \leq C{(1+|\xi|+|\lambda|)}^{-\frac{N}{2}}. \tag2.21$$ 

By a theorem of Greenleaf, (see (0.3) above), 
the inequality $(0.2)$ holds with $q_0=\frac{2(2l+N)}{4l+N}$. 
An application of Lemma 2.3 completes the proof. 
\enddemo 

\demo{Proof of Theorem 2.6} Let $F_0$ be defined as in $(1.7)$. By Lemma 1.4 
$$ |F_{0}(\xi,\lambda)| \leq C{(1+|\xi|+|\lambda|)}^{-\frac{1}{2}}. \tag2.22$$ 

Applying Greenleaf's theorem as above we get 
$q_0=\frac{10}{9}$. An application of Lemma 2.3 completes the proof. 
\enddemo 

\demo{Proof of Theorem 2.9} Let $\Sigma$ be the level set 
$\{x:\phi_2(x)=1\}$, and let
$$
\widehat{d\sigma}(\xi,\,\lambda_1,\,\lambda_2) =
\int_{{\Bbb R}^2} e^{i(
\langle x\,\xi\rangle+\lambda_1\phi_1(x)+\lambda_2\phi_2(x))}
\chi(x)\,dx,
\tag 2.23
$$
where $\chi$ is a smooth cutoff function.
Let  $Z_0$, $Z_1$ and $Z_2$  be defined as in Lemma $2.7$ with respect to 
$\phi_2$. 
Recall that $Z_0\cup Z_1\cup Z_2$ is the union of a finite number of 
lines through the origin.
\par
Let $\{\Gamma_j(x)\}_{j\leq N_1}$ and 
$\{T_j(x)\}_{j\leq N_2}$ be  two finite families of cones in ${\Bbb R}^2$ 
with the following properties:
\item {i)}
Each $\Gamma_j(x)$  contains exactly one line of  
$Z_0\cup Z_1$, and each  $T_j(x)$  contains exactly one line of  
$Z_2$
\item{ii)}
$\Gamma_j \cap \Gamma_i =\{0,\,0\}$ if $i\ne j$,
$T_j \cap T_i =\{0,\,0\}$ if $i\ne j$,
and $\Gamma_j \cap T_i =\{0,\,0\}$.
\par
Let $\alpha_j $ be the characteristic function  of $\Gamma_j$ and let 
$\beta_j$ be the characteristic function of $T_j$.
Then
$$
\widehat{d\sigma}(\xi,\,\lambda_1,\,\lambda_2) =
\int_{{\Bbb R}^2} e^{i(
\langle x\,\xi\rangle+\lambda_1\phi_1(x)+\lambda_2\phi_2(x))}
\chi(x) \left( \sum_{j=1}^{N_1}
\alpha_j(x)\right) \,dx + 
$$
$$
\int_{{\Bbb R}^2} e^{i(
\langle x\,\xi\rangle+\lambda_1\phi_1(x)+\lambda_2\phi_2(x))}
\chi(x) \left(\sum_{j=1}^{N_2}
\beta_j(x)\right) \,dx 
+ 
\tag 2.24
$$
$$
\int_{{\Bbb R}^2} e^{i(
\langle x\,\xi\rangle+\lambda_1\phi_1(x)+\lambda_2\phi_2(x))}\chi(x) 
\left(1-\sum_{j=1}^{N_1}\alpha_j(x)-\sum_{j=1}^{N_2}\beta_j(x)\right)\,dx
= \widehat{d\sigma}_1 + \widehat{d\sigma}_2 + \widehat{d\sigma}_3.
$$
We first consider 
$$
\widehat{d\sigma}_3(\xi,\,\lambda_1,\,\lambda_2) =
\int_{{\Bbb R}^2} e^{i(
\langle x\,\xi\rangle+\lambda_1\phi_1(x)+\lambda_2\phi_2(x))}
\tilde\eta(x)\chi(x)  \,dx,
\tag 2.25
$$
where we have set 
$\tilde\eta(x)= 1-\sum_{j=1}^{N_1}\alpha_j(x)-\sum_{j=1}^{N_2}\beta_j(x)$. On the support of $\tilde\eta$  the curvature of $\Sigma=\{x: \phi_2(x) =1\}$  
never vanishes and $\phi_2$ vanishes only at zero. 
\par
We recall that by the Stein-Tomas observation, (0.2) is equivalent to  the 
inequality
$$
||  \widehat{d\sigma}*f||_{L^{p'}({\Bbb R}^{4})}
\leq C_j||f||_{L^p({\Bbb R}^{4})}.
\tag 2.26
$$ 
Let $\rho$ be a smooth cutoff function supported in the 
interval $(1,\,2)$, such that $\sum_{j=0}^{+\infty} \rho(2^jt) \equiv 1$. 

Let
$$
\widehat{d\sigma_{3,j}}(\xi,\,\lambda_1,\,\lambda_2) =
\int_{{\Bbb R}^2} e^{i(
\langle x\,\xi\rangle+\lambda_1\phi_1(x)+\lambda_2\phi_2(x))}
\rho(2^jx)\tilde\eta(x)  \,dx.
\tag 2.27
$$

The assumptions of Lemma 1.5 are satisfied on the support of 
$\rho(2^jx)\eta(x)$,   hence  the inequality $(2.26)$ holds for the measure 
$d\sigma_{3,j}$ with $p=\frac{6}{5}$.
Since the sharp exponent $p_0$ cannot exceed   $\frac{6}{5}$,  
the estimate  $(2.26)$ holds for $p\leq p_0$, provided that $m_1+m_2 \ge 4$. 

If we make the change of variables sending $x \to 2^{-j}x$,  and if we observe 
that  $\tilde\eta$ is invariant with respect to dilations, we see that
$$ 
\widehat{d\sigma_{3,j}}(\xi,\,\lambda)= 
2^{-2j}\int_{{\Bbb R}^2} e^{i(
\langle 2^{-j}x\,\xi\rangle+\lambda_1 2^{-m_1j}\phi_1(x)+\lambda_2 2^{-m_2j}\phi_2(x))}
\rho(x)\tilde\eta(x)  \,dx = 
$$
$$
2^{-2j}  \widehat{d\sigma_{3,0}}(2^{-j}\xi, \,2^{-m_1j}\lambda_1,\,  
2^{-m_2j}\lambda_2).
\tag 2.28
$$
Without  loss of generality  we can replace $\tilde\eta$ by a function
$\eta\in C ^\infty({\Bbb R}^n\slash\{0\})$,  homogeneous of degree zero, 
whose support coincides with the  support of $\tilde \eta$.

Let $\tau_j$ denote the nonisotropic dilation
$$ \tau_j f(\xi,\,\lambda_1,\,\lambda_2)=
f(2^j\xi,\,2^{m_1j}\lambda_1,\, 2^{m_2j}\lambda_2).\tag 2.29
$$

Then
$$ \widehat{d\sigma_{3,j}}*f= 
2^{-2j}( \tau_{-j}  \widehat{d\sigma_{3,0}}*f)(\xi,\,\lambda_1,\,\lambda_2)=
2^{-2j}( \tau_{-j}  \widehat{d\sigma_{3,0}}*
( \tau_{-j} \tau_jf))(\xi,\,\lambda_1,\lambda_2).\tag 2.30
$$

A change of variables shows that 
$ \tau_{-j}  \widehat{d\sigma_{3,0}}*(\tau_{-j} \tau_jf)=
2^{j(2+m_1+m_2)} \tau_{-j}(\widehat{d\sigma_{3,0}}*\tau_jf).$
It follows that 
$${||{  \widehat{d\sigma_{3,j}}*f}||}_{p'}= 
2^{-2j+ j(2+m_1+m_2)}
{||{ \tau_{-j}(\widehat{d\sigma_{3,0}}*\tau_jf)}||}_{p'}=$$
$$=2^{j(m_1+ m_2)}
2^{\frac{j}{p'} (2+m_1+ m_2)}
{||{  \widehat{d\sigma_{3,0}}*( \tau_jf)}||}_{p'}\leq
\tag 2.31
$$
$$\leq C 2^{j(m_1+ m_2)}
2^{\frac{j}{p'} (2+m_1+m_2)}
{||{ \tau_jf}||}_{p}= 
C2^{j(m_1+ m_2)}
2^{j\left(\frac{1}{ p'}-\frac{1}{p}\right) (2+m_1+ m_2)}
{||{f}||}_{p}.
$$

The series $\dsize
\sum_{j=1}^{\infty}{||{\widehat{d\sigma_{3,j}}*f}||}_{p'}$ 
converges, provided that 
$$m_1+ m_2 +
\left(\frac{1}{p'}-\frac{1}{p}\right) (n+m_1+ m_2)<0,\tag 2.32
$$
which  yields
$\dsize p<p_0$. Hence the measure $d\sigma_3$ satisfies 
the inequality $(2.26)$ with $p < p_0$.  
\bigskip
We consider now   $\dsize
\widehat{d\sigma}_2(\xi,\lambda_1,\,\lambda_2) =\int_{{\Bbb R}^2} 
e^{i(\langle x\,\xi\rangle+\lambda_1\phi_1(x)+\lambda_2\phi_2(x))}
\chi(x) \left(\sum_{j=1}^{N_2}
\beta_j(x)\right) \,dx $.
\par
Let    
$$\widehat{d\mu}(\xi,\lambda_1,\,\lambda_2) =
\int_{{\Bbb R}^2} 
e^{i(\langle x\,\xi\rangle+\lambda_1\phi_1(x)+\lambda_2\phi_2(x))}
\chi(x)\beta_j(x) \,dx,
\tag 2.33
$$ 
where   $j\leq N_2$ is fixed.
\par
In the proof of lemma 1.4 we observed that  
${{\phi_1(x)}_{|}}_{\Sigma} = c $  implies that
$\phi_1(x)=c\phi_2(x)^{\frac{m_1}{m_2}}$ for every $x\in {\Bbb R}^2$.
Observe that $\Sigma$ is star-shaped with respect to the origin, because 
if $x_0\in \Sigma$, then for every $t\leq1$, $\phi_1(tx_0)=t^{m_1}\leq 1$. 
In the polar coordinates associated to $\Sigma$,

$$
\widehat{d\mu}(\xi,\,\lambda_1,\,\lambda_2) =
\int_{\Sigma}\int_0^{+\infty}
e^{i(r\langle \omega\,\xi\rangle+c\lambda_1 r^{m_1}   +\lambda_2 r^{m_2}  )}
r\chi(r\omega)\beta_j(\omega)\,dr\,d\omega.
\tag 2.34
$$
Without loss generality  $\chi$ is radial. Consider 
the (unique) point of 
$Z_2\cap \Sigma\cap T_j$, which can be taken to be $(0,\,0)$. Suppose that 
$\Sigma\cap T_j$ is supported in a sufficiently small neighborhood  of  
$(0,\,0)$.  Since $\Sigma$ is finite type $m$, it 
can be   written  as the graph of a smooth 
function $\psi(t)=t^mg(t)$, where $g(0) \ne 0$, and 
$$
\widehat{d\mu} (\xi,\,\lambda_1,\,\lambda_2) =
\int_{{\Bbb R}}\int_0^{+\infty}
e^{i(rt\xi_1 + rt^mg(t) \xi_2 +c\lambda_1 r^{m_1}  +\lambda_2 r^{m_2})}
r\chi(r)\beta_j(t,t^mg(t))\,dr\,dt.
\tag 2.35
$$
Let $\rho\in C^{\infty}({\Bbb R}^2)$ be supported in 
$(\frac{1}{4},\, 4)$, such that   $\rho \equiv 1$ in $(1,\,2)$, and   
$\sum_{j=1}^{+\infty}\rho(2^jt)\equiv 1$.  Let 
$$
\widehat{d\mu_j }(\xi,\,\lambda_1,\,\lambda_2)  =
\int_{{\Bbb R}}\int_0^{+\infty}
e^{i(rt\xi_1 + rt^mg(t)\xi_2 +c\lambda_1 r^{m_1}  +\lambda_2 r^{m_2})}
 \,\rho(2^jt)\chi(r)rdr\,dt.\tag 2.36
$$ 
The integral with respect to $t$  is supported over a 
dyadic piece of $\Sigma$ where the  Gaussian curvature does not vanish. 
By Lemma 1.5 and the Stein-Thomas observation, the estimate $(2.26)$
holds for the measure $d\mu_j$, for $p\leq p_0$, with a constant $C_j$.
\par
In order to estimate $C_j$   we make a change of variables  in the expression
for $\widehat{ d\mu_j}$  setting $s=2^jt$. We have

$$
\widehat{d\mu_j }(\xi,\,\lambda_1,\,\lambda_2)  =
2^{-j}\int_{{\Bbb R}}\int_0^{+\infty}
e^{i(r2^{-j}s\xi_1 + r2^{-mj}s^mg(2^{-j}t)\xi_2 +
c\lambda_1 r^{m_1}  +\lambda_2 r^{m_2})}
\chi(r)\,\rho(s)rdr\,ds.
\tag 2.37
$$ 
Let $\tau_j$ be the nonisotropic dilation
\par\noindent
$\tau_j(f)(x_1,\,x_2,\,x_3,\,x_4) 
= f(2^{-j}x_1, 2^{-mj}\,x_2,\,x_3,\,x_4)$ and let 
$$
\widehat{d\widetilde\mu_j }  =
\int_{{\Bbb R}}\int_0^{+\infty}
e^{i(rs\xi_1 + rs^mg(2^{-j}t)\xi_2 +c\lambda_1 r^{m_1}  +\lambda_2 r^{m_2})}
\chi(r)\,\rho(s)rdr\,ds.
\tag 2.38
$$ 
Then 
$$
\widehat{d\mu_j}*f = 2^{-j}(\tau_{j}(\widehat{d\widetilde\mu_j})*f)
=2^{-j}(\tau_{j}(\widehat{d\widetilde\mu_j})*\tau_{j}(\tau_{-j}f)).
\tag 2.39
$$
A change of variables shows that 
$\tau_{j}(\widehat{d\widetilde\mu_j})*\tau_{j}\tau_{-j}f$=
$ 2^{j(m+1)}\tau_{j}(\widehat{d\widetilde\mu_j}* \tau_{-j}f)$, and that
\newline $||\tau_j\psi||_{L^{q}({\Bbb R}^4)} =
2^{j\frac{(m+1)}{q}} ||\psi||_{L^{q}({\Bbb R}^4)}$ for every 
$\psi \in L^{q}({\Bbb R}^4)$ and $q\ge 1$. 
Then we can write the following string of inequalities:
$$
||\widehat{d\mu_j}*f||_{L^{p\prime}({\Bbb R}^4)} = 
2^{jm}||\tau_{j}(\widehat{d\widetilde\mu_j}* \tau_{-j}f)||_{L^{p\prime}
({\Bbb R}^4)}
=$$
$$
2^{jm +j\frac{(m+1)}{p\prime}}||  
\widehat{d\tilde\mu_j}* \tau_{j}f||_{L^{p'}({\Bbb R}^4)}
\leq 2^{jm +j\frac{(m+1)}{p\prime}}C_j||\tau_{j}f||_{L^{p}({\Bbb R}^4)} =
\tag 2.40
$$
$$
2^{jm +j(m+1)\left(\frac{1}{p\prime}-\frac{1}{p}\right)}
C_j||f||_{L^{p}({\Bbb R}^4)}.
$$
We must prove that the constants $C_j$ in the above expression are uniformly 
bounded. In fact the sum 
$\dsize 
\sum_{j=1}^{+\infty}||\widehat{d\mu_j}*f||_{L^{p'}({\Bbb R}^4)}$
converges if 
$ m+(m+1)\left(\frac{1}{p'}-\frac{1}{p}\right)<0$, hence if 
$ p< \frac{2(m+1)}{2m+1}$. 
Since $\frac{2(m+1)}{2m+1} \ge p_0$ when $m_1+m_2\ge 2m$, the 
estimate $(2.26)$ holds for the measure $d\mu$, and 
consequently for the measure $d\sigma_2$,  
for $p\leq p_0$, provided  that $m_1+m_2 \ge 2m$, and with 
$p=\frac{2(m+1)}{2m+1}$ if $m_1+m_2<2m$.

>From the proof of the theorem of Greenleaf it follows that the bounds for the 
constants $C_j$  depend on a finite number of  derivatives of the phase 
function of $\widehat{d\tilde\mu_j}$,  $\Phi_j(r,\,s)=
rs\xi_1 + rs^mg(2^{-j}t)\xi_2 +c\lambda_1 r^{m_1}  +\lambda_2 r^{m_2}$.
Since $\Phi_j$ is a smooth function, then, for $j$ large,  
$D^\beta\Phi_j(r,\,s) \approx
D^\beta(rs\xi_1 + rs^mg(0)\xi_2 +c\lambda_1 r^{m_1}  +\lambda_2 r^{m_2})$.
This shows that the constant $C_j$'s are uniformly 
bounded.
\bigskip
We now consider  $\dsize
\widehat{d\sigma}_1(\xi,\lambda_1,\,\lambda_2) =\int_{{\Bbb R}^2} 
e^{i(\langle x\,\xi\rangle+\lambda_1\phi_1(x)+\lambda_2\phi_2(x))}
\chi(x) \left(\sum_{j=1}^{N_1}
\alpha_j(x)\right) \,dx $.
\par
Fix  $j\leq N_1$. After perhaps a 
rotation of coordinates we may assume that  
${{\phi_2}_|}_{\Gamma_j}$ vanishes along the $x_2$ axis.
Then 
$\phi_2(x_1,\, x_2)$ can be  written as $ x_1^Mg(x_1,x_2)$, where 
$g$  does not vanish on the $x_2$ axis, (except perhaps at the origin), 
if  $\Gamma_j\cap S^1$ is  small enough.
\par
Let
$$\widehat{d\mu}(\xi,\,\lambda_1,\,\lambda_2)=\int_{{\Bbb R}^2} e^{i(
\langle x\,\xi\rangle+c\lambda_1x_1^{\gamma M}g(x)^\gamma+
\lambda_2x_1^{M}g(x))}
\chi(x) \alpha_j(x) \,dx,
\tag 2.41
$$
where we have set $\gamma =\frac{m_1}{m_2}$.
Let $\rho\in C^{\infty}({\Bbb R}^2)$ 
be a cutoff  function  supported in $(\frac{1}{4}, \, 4)$, 
such that $\rho\equiv 1 $ in $(1,\, 2)$ and  
$\sum_{j=1}^{+\infty}\rho(t)\equiv 1$. Let 
$$
\widehat{d\mu_j}(\xi,\lambda)=\int_{{\Bbb R}^2} e^{i(
\langle x\,\xi\rangle+
c\lambda_1x_1^{\gamma M}g^\gamma(x)+\lambda_2x_1^{M }g(x))}
\rho\left(2^{j}\frac{x_1}{x_2}\right) \chi(x)\,dx.
\tag 2.42
$$
The above integral   
is  defined over a cone of ${\Bbb R}^2$  where the curvature of $\Sigma$
never vanishes, and $\phi_2$ vanishes only at the origin.  
\par
Let $\rho_0\in C^{\infty}({\Bbb R}^2)$  be  supported in  
$(\frac{1}{4},\,4)$, such that  $\rho \equiv 1$ in $(1,\,2)$,
and $\sum_{k =0}^{+\infty} \rho_0(2^kt) \equiv 1$.

Let
$$
\widehat{d\mu_{j,k}}(\xi,\lambda)=\int_{{\Bbb R}^2} e^{i(
\langle x\,\xi\rangle+
c\lambda_1x_1^{\gamma M}g^\gamma(x)+\lambda_2x_1^{M }g(x))}
\rho\left(2^{j}\frac{x_1}{x_2}\right) \rho_0(2^kx)\,dx.
\tag 2.43
$$
The assumptions of Lemma 1.5 are satisfied on the support of 
$\rho_0(2^kx)\rho\left( 2^j \frac{x_1}{x_2}\right)$,  and hence 
the estimate $(2.26)$ holds for 
the measure $d\mu_{j,k}$ with $p\leq p_0$. 

If we make the change of variables sending $x \to 2^{-k}x$, 
$$ 
\widehat{d\mu_{j,k}}(\xi,\,\lambda)= 
2^{-2k}\int_{{\Bbb R}^2} e^{i(
\langle 2^{-k}x\,\xi\rangle+\lambda_1 2^{-m_1k}\phi_1(x)+\lambda_2 2^{-m_2k}\phi_2(x))}
\rho\left(2^{j}\frac{x_1}{x_2}\right) \rho_0(x)\,dx = 
$$
$$
2^{-2k}  \widehat{d\mu_{j,0}}(2^{-k}\xi, \,2^{-m_1k}\lambda_1,\,  
2^{-m_2k}\lambda_2).
\tag 2.44$$

Let $\tau_k$ denote the nonisotropic dilation
$$ \tau_k f(\xi,\,\lambda_1,\,\lambda_2)=
f(2^k\xi,\,2^{m_1k}\lambda_1,\, 2^{m_2k}\lambda_2).\tag 2.45
$$

Then
$$ \widehat{d\mu_{j,k}}*f= 
2^{-2k}( \tau_{-k}  \widehat{d\mu_{0,k}}*f)(\xi,\,\lambda_1,\,\lambda_2)=
2^{-2k}( \tau_{-k}  \widehat{d\mu_{0,k}}*
( \tau_{-k} \tau_kf))(\xi,\,\lambda_1,\lambda_2).\tag 2.46
$$

A change of variables shows that 
$ \tau_{-k}  \widehat{d\mu_{j,k}}*(\tau_{-k} \tau_kf)=
2^{k(2+m_1+m_2)} \tau_{-k}(\widehat{d\mu_{j,k}}*\tau_kf).$
It follows that 
$${||{  \widehat{d\mu_{j,k}}*f}||}_{p'}= 
2^{-2k+ k(2+m_1+m_2)}
{||{ \tau_{- k}(\widehat{d\mu_{j,k}}*\tau_kf)}||}_{p'}=$$
$$=2^{k(m_1+ m_2)}
2^{\frac{k}{p'} (2+m_1+ m_2)}
{||{  \widehat{d\mu_{j,k}}*( \tau_kf)}||}_{p'}\leq
\tag 2.47
$$
$$\leq C 2^{k(m_1+ m_2)}
2^{\frac{k}{p'} (2+m_1+m_2)}
{||{ \tau_ kf}||}_{p}= 
C2^{k(m_1+ m_2)}
2^{k\left(\frac{1}{ p'}-\frac{1}{p}\right) (2+m_1+ m_2)}
{||{f}||}_{p}.$$

The series $\dsize
\sum_{k=1}^{\infty}{||{\widehat{d\mu_{j,k}}*f}||}_{p'}$ converges,  
provided that 
$$m_1+ m_2 +
\left(\frac{1}{p'}-\frac{1}{p}\right) (n+m_1+ m_2)<0,
\tag 2.48
$$
which  yields
$\dsize p<p_0$. 
The above argument shows  that we can assume that the measure $d\mu_j$ 
is supported away from zero. Hence,
by Lemma 1.5  and  the above observation, the inequality $(2.26)$ 
holds for the measure $d\mu_j$, for   $p\leq p_0$, 
with a constant $C_j$.
In order to estimate $C_j$ 
we perform the change of variables in  $(2.42)$ sending 
$x_1\to 2^{-j}x_1$, $x_2\to x_2$. We  obtain 
$$
\widehat{d\mu_j}(\xi_1,\xi_2,\lambda_1,\lambda_2)= 2^{-j}
\widehat{d\widetilde \mu_j}( 2^{-j}\xi_1, \xi_2,\,
2^{-M\gamma j}\lambda_1,\, 2^{-Mj }\lambda_2),
\tag 2.49
$$
where we have set 
$$
\widehat{d\widetilde\mu_j}
(\xi,\lambda)=
\int_{{\Bbb R}^2} e^{i(
\langle x\,\xi\rangle+c\lambda_1x_1^{\gamma M }
g^\gamma(2^{-j}x_1, x_2)+\lambda_2x_1^{M}g (2^{-j}x_1, x_2))}
\rho\left(\frac{x_1}{x_2}\right)\chi(2^{-j}x_1, x_2)  \,dx.
\tag 2.50
$$
Let $\tau_{j}$ be the nonisotropic dilation 
$\tau_{j}f(x_1,\,x_2,\,x_3,\,x_4)$ 
$ = f(2^{-j}x_1, \,x_2,\,2^{-M\gamma j}x_3,2^{-M j}\,x_4)$.
Then 
$$
\widehat{d\mu_j}*f = 2^{-j}(\tau_{j}(\widehat{d\widetilde\mu_j})*f)
  = 2^{-j}(\tau_{j}(\widehat{d\widetilde\mu_j})*\tau_{j}(\tau_{-j}f)).
\tag 2.51
$$
A change of variables shows that 
$$
\tau_{j}(\widehat{d\widetilde\mu_j})*\tau_{j}(\tau_{-j}f)  =
 2^{j(1+M(1+\gamma))}
\tau_{j}(\widehat{d\widetilde\mu_j}* \tau_{-j}f),
\tag 2.52
$$
and that
$||\tau_{j}g||_{L^{q}({\Bbb R}^4)} =
2^{j\frac{(1+M(1+\gamma))}{q}} ||g||_{L^{q}({\Bbb R}^4)}$ for every $q\ge 1$. 
Then we can write the following string of inequalities:
$$|| \widehat{d\mu}_j *f||_{L^{p\prime}({\Bbb R}^4)} = 2^{j(M(1+\gamma))}||\tau_{j}(\widehat{d\widetilde\mu_j}* \tau_{-j}f)||_{L^{p\prime}({\Bbb R}^4)}
=$$
$$
2^{jM(1+\gamma) +j\frac{(1+ M(1+\gamma))}{p\prime}}
|| \widehat{d\tilde\mu_j}* \tau_{j}f||_{L^{p\prime}({\Bbb R}^4)}\leq
2^{j(M(1+\gamma)) +j\frac{M(1+\gamma)+1}{p\prime} }
C_{j}||\tau_{j}f||_{L^{p}({\Bbb R}^4)} =
\tag 2.53
$$
$$
2^{jM(1+\gamma) +j(M(1+\gamma)+1)\left(\frac{1}{p\prime}-\frac{1}{p}\right)}
C_{j}||f||_{L^{p}({\Bbb R}^4)}.
$$
If we show that the constants $C_{j}$ in the above expression are uniformly 
bounded then we are done. In fact the sum $\dsize 
\sum_{j=1}^{+\infty}||\widehat{d\mu_j}*f||_{L^{p'}({\Bbb R}^4)}$
converges if 
\par\noindent
$
M(1+\gamma) +(M(1+\gamma)+1)\left(\frac{1}{p\prime}-\frac{1}{p}\right) 
<0$
hence if
$$
p <\frac{2(1+M(\gamma+1))}{1+2M(\gamma+1)}.
\tag 2.54
$$ 
Since $\frac{2(1+M(\gamma+1))}{1+2M(\gamma+1)}\ge p_0$  when  
$m_1+m_2\ge 2M(1+\gamma)$, 
the estimate $(2.26)$ holds for the 
measure $d\mu$, and consequently for the  
measure $d\sigma_1$,   with   $p< p_0$, provided 
that $m_1+m_2 \ge 2M(1+\gamma)$, and with  
$ p < \frac{2(1+M(\gamma+1)}{1+2M(\gamma+1)}$ if $4< m_1+m_2<2M(1+\gamma)$.  
\par
By a theorem of Greenleaf the bounds for the constants
$C_j$ in  $(2.53)$  depend only on a finite number of  derivatives of 
the phase function of $ \widehat{d\tilde\mu_j}$,  
$  \Phi_j(x)= (c\lambda_1x_1^{\gamma M} g^\gamma(2^{-j}x_1, x_2) +
\lambda_ 2x_1^{ M} g(2^{-j}x_1, x_2)) $
Since  the above function  is  smooth, then, for $j$ large, 
$ D^{\beta} \Phi_j(x)\approx D^{\beta}(c\lambda_1x_1^{\gamma M} 
g^\gamma(0, x_2) +
\lambda_ 2x_1^{M} g(0, x_2))$. 
This shows that the constants $C_j$ are uniformly 
bounded, thus concluding the proof of the theorem.  
\enddemo

\vskip.25in

\head Section 3 \endhead 
\head Restriction theorems- Non-isotropic estimates \endhead
\vskip.125in

\proclaim{Theorem 3.1} Let $S=\{(x,x_{n+1},...,x_{n+l}) \in {\Bbb R}^{n+l}: 
x_{n+1}={\phi}_1(x),..., x_{n+l}={\phi}_l(x)\}$, where $l=\frac{n(n+1)}{2}$, and
the ${\phi}_j$ denote the distinct monomials of degree $2$. Then the estimate
$(0.2)$ holds with the sharp exponent $p_0=\frac{2(n+2)}{2n+3}$. 
\endproclaim

\proclaim{Theorem 3.2} Let $S$ denote a compact piece of the manifold 
$\{(x,x_{n+1},...,x_{n+l}) \in {\Bbb R}^{n+l}: 
x_{n+1}={\phi}_1(x),..., x_{n+l}={\phi}_l(x)\}$,where $\phi_j \in C^{\infty}
({\Bbb R}^n\slash 0)$
is homogeneous of degree $m \ge 2n$. Suppose that no linear combination 
of the  $\phi_j$'s  vanishes on a subset of postive measure of $S^{n-1}$.
Let ${\Phi}(x)= ({\phi}_1(x),..., {\phi}_l(x))$. 
Suppose that $\Phi(\omega)\not=(0,...,0)$, 
$\omega \in S^{n-1}$. 
Then the estimate $(0.2)$ holds with the 
sharp exponent $p_0$ given by $(0.4)$. 
\endproclaim  

\remark{Remark 1} The restriction $m\ge 2n$ in Theorem 3.2 is not necessary. 
In fact, using the techniques in (\cite{IS}), one can prove Theorem 3.2 under the  
weaker restriction $m\ge n$. 
\endremark 

\remark{Remark 2} Theorem 3.2 implies the natural generalization of Theorem 3.1
to the case where $l=C^{m+n-1}_{m}$, $(C^a_b=\frac{a!}{b!(a-b)!})$, and the 
${\phi}_j$ are the disitinct monomials of degree $m$. 
\endremark 
\vskip.25in

\demo{Proof of Theorem 3.1} 
Let $\lambda = (\lambda_1,\cdots, \lambda_l)$, and let
$$
\widehat{d\sigma}(\xi,\,\lambda) =
\int_{{\Bbb R}^n} e^{i(
\langle x\,\xi\rangle+ \lambda_1{\phi}_1+\cdots \lambda_l{\phi}_l(x))}dx.
$$
Let $A_\lambda$ be the matrix associated to the quadratic form  
$ \lambda_1{\phi}_1(x)+\cdots \lambda_l{\phi}_l(x)$, and 
$$
\widehat{d\sigma}(\xi,\,\lambda) =
\int_{{\Bbb R}^n} e^{i(
\langle x\,\xi\rangle+ \langle x, \, A_\lambda x\rangle)}dx.
\tag 3.1
$$
Thus, $\widehat{d\sigma}(\xi,\,\lambda)$ is the Fourier transform of 
$\dsize e^{i \langle x, \, A_\lambda x\rangle}$,  and  an easy generalization 
of the well-known formula  for the  Fourier transform of the Gaussian
functions, (see e.g.\cite{WR} pg. $186$), yields
$$
\widehat{d\sigma}(\xi,\,\lambda) = 
\frac{(2\pi)^{\frac{n}{2}}}{|\text{det}A_\lambda|^{1/2}}
e^{-\frac{i}{2}\langle \xi,\, A_\lambda^{-1} \xi\rangle+\frac{\pi i}{4}
\text{sign}(A_\lambda)}.
\tag 3.2
$$
Let
$$
K_z(\xi,\lambda)=\psi(z)
\text{det}(A_\lambda)^{z}\widehat{d\sigma}(\xi,\,\lambda)
=
\psi(z)(2\pi)^{\frac{n}{2}}|\text{det}A_\lambda|^{z-1/2}
e^{-\frac{i}{2}\langle \xi,\, A_\lambda^{-1} \xi\rangle +\frac{\pi i}{4}
\text{sign}(A_\lambda)},
\tag 3.3
$$
with 
$$
\psi(z)=(2\pi)^{nz+\frac{n(n+1)}{2}} \Pi_{j=0}^{n-1}\
\Gamma^{-1}(z+1+\frac{j}{2}) 
e^{-i\pi\frac{n}{2}(z+\frac{n+1}{2})},
\tag 3.4
$$ 
where $\Gamma$ is the standard Gamma function.
Let $T_z(f)(\xi,\,\lambda)= (f*K_z)(\xi,\,\lambda)$.
\par
We will prove that $T_z$ is an continuous family of operators when 
Re$(z) \in [-\frac{n+1}{2},\,  \frac{1}{2}]$, is analytic  
when Re$(z) \in (-\frac{n+1}{2}\, \frac{1}{2})$,
and that
\bigskip
\item{i)}
$||T_z(f)||_{L^\infty({\Bbb R}^n)} \leq C_1(z)||f||_{L^1({\Bbb R}^n)} $, 
when Re$(z) =\frac{1}{2}$,
\item{ii)}
$||T_z(f)||_{L^2({\Bbb R}^n)} \leq C_2(z)||f||_{L^2({\Bbb R}^n)}$, when 
Re$(z)=- \frac{n+1}{2}$,
 \item{iii)}
$C_1(z)$ and $C_2(z)$ have at most exponential growth with respect to Im$(z)$. 
\bigskip
\noindent
Stein's analytic
interpolation theorem (see e.g. \cite{So})  will then imply  that 
\par
\noindent
$T_z : L^{p_0}({\Bbb R}^{n+l})\to L^{p_0\prime}({\Bbb R}^{n+l})$ 
is a bounded operator  
when  $z=0$. 
Since $T_z(f)= f*\widehat{d\sigma}$, the Stein-Thomas observation 
(see the proof of Theorem 2.9) implies the conclusion of our theorem.
\bigskip
To prove $(i)$  
we observe that when $z =\frac{1}{2}+iy$, with $y\in {\Bbb R}$, then   
\newline $K_z(\xi,\lambda)=
\widehat{d\sigma}(\xi,\,\lambda)\text{det}(A_\lambda)^{\frac 1 2 +iy}
\psi(\frac{1}{2}+iy)$  is bounded by $C_1(y)=\pi^{\frac{n}{2}} 
|\psi(\frac{1}{2}+iy)|$.
By the Hausdorff -Young inequality we have that $||f*K_{1+iy}||_{\infty}\leq C_1(y)||f||_{1}$. $ (i)$  is then  satisfied, and one can check, using 
Stirling's formula, that $C_1(y)$ has at most exponential growth.
\medskip
To prove $(ii)$ it is enough to show that  
$|\widehat K_z|$ is a bounded function when Re$(z)=-\frac{n+1}{2}$.
To compute the the Fourier transform of 
$K_z(\xi,\,\lambda)$ with respect to $\xi $, 
${\Cal F}_\xi( K_z)(x,\,\lambda)$,  we use again the formula $(3.3)$ obtaining
$$
{\Cal F}_\xi( K_z)(x,\,\lambda)= \psi(z)
|\text{det}A_\lambda|^{z}e^{i\langle x,\, A_{\lambda}x\rangle}.
\tag 3.5
$$
Hence 
$$
\widehat{K_z}(x, x') = \psi(z) 
\int_{{\Bbb R}^l} e^{i(
\langle x',\,\lambda\rangle+ \langle x, \, A_\lambda x\rangle)}
\text{det}(A_\lambda)^zd\lambda,
\tag 3.6
$$
where we have set $x'=(x_{n+1},\cdots x_{n+l})$. We recall that the above 
identities hold in distribution sense.
\par
Since the phase  of the above integral is  a linear  function of $\lambda$,
we reduce to computing the Fourier transform of 
det$(A_\lambda)^z$. 
We need the following lemma (see \cite{TS},  pg. $48$). 

\proclaim{Lemma 3.3} Let  $V_{\Bbb R}$  be the space of the
real and symmetric matrices and let $V_i\subset V_{\Bbb R}$ 
be the subset of  the   matrices with $i$ positive and $n-i$ negative eigenvalues. Let $G_i(z)$ be the distribution 
$$
G_i(z)(f) = \int_{V_i}f(Y)|\text{det}(Y)|^zdY
\tag 3.7
$$
where $dY$ is the standard Euclidean measure on $V_{\Bbb R}$. 
Then the distribution  $G_i(z)$, viewed as a function of $z$, 
has analytic continuation to a meromorphic function in the whole complex plane 
satisfying 
$$ 
\widehat{G_i(z)}
= 
\psi^{-1}(z)
\sum_{j=0}^{n} c_{i,j}(z)G_j\left(-z-\frac{n+1}{2}\right),
\tag 3.8
$$
where $\psi$ is as in $(3.4)$ and the $C_{i,j}(z)$ are bounded coefficients.
\endproclaim
>From the above formula we deduce that, modulo bounded constants, 
$$
{\Cal F}(\text{det}(A_\lambda)^z )(\eta)= 
\psi^{-1}(z)\text{det}(A_{\eta})^{-z-\frac{n+1}{2}},
\tag 3.9
$$
where the above formula holds in distribution sense.
Since Re$(z) =\frac{n+1}{2}$, the above is a bounded function of $\eta$.
This shows that $\widehat{K_z}$ is a bounded function of $(x,\, x')$, and 
completes the proof  of the theorem.
\enddemo 

\demo{Proof of Theorem 3.2} Let
$$
\widehat{d\sigma}(\xi,\,\lambda) =
\int_{{\Bbb R}^n} e^{i(
\langle x\,\xi\rangle+ \lambda_1\phi_1(x)+\lambda_l\phi_l(x))}
\chi(x)\,dx,
\tag 3.10
$$
where $\chi$ is a smooth cutoff function.
By Theorem $1.1B$  we can write
$$
\widehat{d\sigma}(\xi,\,\lambda) = \widetilde{B}(\xi,\lambda)
\tag 3.11
$$
where $\dsize |\widetilde{B}(\xi,\lambda)|\leq \int_{S^{n-1}} \frac
{d\omega}{|\sum_{j=1}^{l} \phi_j(\omega)\lambda_j|^{\frac n m}}
$. 

Let 
$$
\dsize B(\lambda) =\int_{S^{n-1}} \frac
{d\omega}{|\sum_{j=1}^{l}\phi_j(\omega)\lambda_j|^{\frac n m}}, 
\tag 3.12
$$  
$$
K_z(\xi,\lambda)=
\widehat{d\sigma}(\xi,\,\lambda)  {B}^{z}(\lambda)
\psi(z),
$$
with 
$$
\psi(z) = \Gamma^{-1}(-\frac{zn}{2m}+\frac{1}{2}) 
\Gamma^{-1}(-\frac{nz}{m}+l),
$$
where $\Gamma$ is the standard Gamma function.
Let $T_z(f)(\xi,\,\lambda)= (f*K_z)(\xi,\,\lambda)$.

We will prove that $T_z$ is a continous family of operators when 
Re$(z) \in [-1,\,\frac{lm}{n}]$, is analytic  
when Re$(z) \in (-1,\,\frac{lm}{n})$,
and that
\bigskip
\item{i)}
$||T_z(f)||_{L^\infty({\Bbb R}^n)} \leq C_1(z)||f||_{L^1({\Bbb R}^n)} $, 
when Re$(z) = -1$,
\item{ii)}
$||T_z(f)||_{L^2({\Bbb R}^n)} \leq C_2(z)||f||_{L^2({\Bbb R}^n)}$, when 
Re$(z)=\frac{ml}{n}$,
 \item{iii)}
$C_1(z)$ and $C_2(z)$ have at most exponential growth with respect to Im$(z)$. 
\bigskip
\noindent
Stein's analytic
interpolation theorem (see e.g. \cite{So})  will then imply  that 
\par
\noindent
$T_z : L^{p_0}({\Bbb R}^{n+l})\to L^{p_0\prime}({\Bbb R}^{n+l})$ 
is a bounded operator  when  $z=0$. 
Since $T_z(f)= f*\widehat{d\sigma}$, Stein-Thomas observation 
(see the proof of Theorem 2.9) implies the conclusion of our theorem.
\bigskip
To prove the estimate $(i)$  
we observe that when $z =-1+ iy$, with $y\in {\Bbb R}$, then   
\newline $K_z(\xi,\lambda)=
\widehat{d\sigma}(\xi,\,\lambda)   
{B}^{iy -1}(\lambda) \psi(-1+iy)$  
is bounded by 
$C_1(y)= C|\psi_1(-1+iy)|$.
By the Haussdorf -Young inequality we have that 
$||f*K_{1+iy}||_{\infty}\leq C_1(y)||f||_{1}$. 
The estimate $ (i)$  is then  satisfied, and one can 
check using Stirling's formula that $C_1(y)$ has at most exponential growth.
\medskip
To prove the estimate $(ii)$ it is enough to show that $\widehat K_z$ is a
bounded function  when Re$(z)=\frac{ml}{n}$. 
Since $d\sigma$ is a finite measure 
it is enough to prove that $\psi(z)\widehat{{B}^{z}}(\lambda)$ is bounded.
\par
Let $X(\omega)$ be the vector defined by the equation $\dsize 
\sum_{j=1}^{l}\lambda_j\phi_j(\omega) =\langle \lambda,\, X(\omega)\rangle$.
Then
$$
\widehat{B_z}(y)=\psi(z)\int_{ {\Bbb R}^l} \left(
\int_{S^{n-1}}\frac{d\omega}{|\langle \lambda,\, X(\omega)\rangle|^
{\frac{n}{m}}}
\right)^{z} e^{i\langle \lambda,\, y \rangle} d\lambda.
\tag 3.13
$$
Since $\widehat{B_z}(y)$ is homogeneous of degree zero with respect to $y$,
we can assume  $|y|=1$.

In polar coordinates with respect to $\lambda$,
with $\lambda =r\eta$, we have
$$
\widehat{B_{z}}(y)=\psi(z)\int_{0}^{+\infty}
r^{-\frac{n}{m}z+l-1}
\int_{S^{l-1}}\left(
\int_{S^{n-1}}
\frac{ e^{\frac{i}{z}r\langle \eta,\, y \rangle}d\omega}
{|\langle \eta,\,X(\omega)\rangle|^{\frac{n}{m}}}
\right)^{z} d\eta dr.
\tag 3.14
$$
Let $\chi(r)\in {\Cal C}^{\infty}( {\Bbb R})$ be such that 
$\chi(r) \equiv 1$, when $r\in (0,\, 1)$, and 
$\chi(r) \equiv 0$ when $r \in (2,\, +\infty)$. Let
$$
\widehat{B_{z,M}}(y)= 
\psi(z)\int_{0}^{+\infty}
r^{-\frac{n}{m}z+l-1}\chi\left(\frac{r}{M}\right)
\int_{S^{l-1}}\left(
\int_{S^{n-1}}
\frac{ e^{\frac{i}{z}r\langle \eta,\, y \rangle}d\omega}
{|\langle \eta,\,X(\omega)\rangle|^{\frac{n}{m}}}
\right)^{z} d\eta dr.
\tag 3.15
$$
We will prove that  $\widehat{B_{z,M}}$ is bounded by a constant $C$
independent of $M$ and that 
$\widehat{B_{z,M}} \to \widehat{B_{z}}$ in distribution sense as 
$M\to \infty$.
>From the above it  follows that $\widehat{B_z}$  is bounded. 
Indeed, since the balls are sequentially compact in the 
weak $*$ topology  in $L^{\infty}( {\Bbb R}^l)$, there exists a sequence
$\{\widehat{B_{z,M_j}}\}_{j\in {\Bbb N}}\subset  
\{\widehat{B_{z,M}}\}$  which converges to a bounded function  in the
weak $*$ topology of $L^{\infty}$, and hence converges also in 
distribution sense.
Consequently, $\widehat{B_z}= \lim_{j\to +\infty}\widehat{B_{z, M_j}}$, 
which is a bounded function.

Recalling that  by  assumption 
the vector $X(\omega)$ is never  zero on $ S^{n-1}$, 
we can construct an orthogonal matrix $A_\omega$ with the 
property that $ A_\omega \frac{X(\omega)}{|X(\omega)|} = 
(1, \, 0, \, \cdots,\, 0)$.
The first row of $A_\omega$ is $\frac{X(\omega)}{|X(\omega)|} $ and the 
other rows  are a set of $l-1$ vectors which, 
togeter with $\frac{X(\omega)}{|X(\omega)|}$, 
determine an orthonormal  basis of ${\Bbb R}^l$ for every $\omega \in S^{n-1}$.

We make the change of variables 
$\eta\to A_\omega^t\eta$ in the expression for $\widehat{B_{M,z}}$.
Since $A_\omega$ is orthogonal, 
the change of variables maps $S^{l-1}$ into itself and the determinant of the 
Jacobian matrix of the transformation is $1$. We obtain
$$
\widehat{B_{z,M}}(y)=\psi(z)\int_{0}^{+\infty}
r^{-\frac{n}{m}z+l-1}\chi\left(\frac{r}{M}\right) 
\!\int_{S^{l-1}}\!\left(
\int_{S^{n-1}}\!\frac{|X(\omega)|^{\frac{n}{m}}}
{|\langle \eta,\,A_\omega X(\omega)\rangle|^{\frac{n}{m}}}
e^{\frac{i}{z}r\langle \eta,\, A_\omega y \rangle}d\omega \!\right)^{z} d\eta dr
$$
$$
=\psi(z)\int_{0}^{+\infty}
r^{-\frac{n}{m}z+l-1} \chi\left(\frac{r}{M}\right)
\int_{S^{l-1}} \frac{1}{|\eta_1|^{\frac{nz}{m}}}\left(
 \int_{S^{n-1}}{|X(\omega)|^{\frac{n}{m}}}
e^{\frac{i}{z}r\langle \eta,\, A_\omega y \rangle} d\omega\right)^{z}d\eta dr. 
\tag 3.16
$$
The integral with respect to $\omega$  is a continous function of  $r$, $\eta$, 
$y$, say $F(r,\,\eta,\, y)$. 
In particular $F(0,\,\eta,\,y)=  
\int_{S^{n-1}}{|X(\omega)|^{\frac{n}{m}}}d\omega<\infty.$ 

Since Re$(z)\ge l$ by assumption, we can check that we can compute at least
$l-1$ derivatives of  $F(r,\,\eta,\,y)^z$ with respect to $r$ and $\eta$. Then
$$
\widehat{B_{z,M}}(y)= \psi(z)\int_{0}^{+\infty}
r^{-\frac{n}{m}z+l-1} \chi\left(\frac{r}{M}\right)
\int_{S^{l-1}} \frac{1}{|\eta_1|^{\frac{nz}{m}}}F^z(r,\,\eta,\,y) drd\eta= 
$$
$$
\psi(z)\int_{0}^{+\infty}r^{-\frac{n}{m}z+l-1}\chi\left(\frac{r}{M}\right)
\int_{-1}^{1}
\int_{\{|\eta^\prime |=(1-|\eta_1|^2)^{\frac{1}{2}}\}} 
\frac{F^z(r,\,\eta,\,y)}{|\eta_1|^{\frac{nz}{m}}} drd\eta_1d\eta^\prime.
\tag 3.17
$$
Here  we set $\eta = (\eta_1,\, \eta^\prime)$, and we let $d\eta^\prime$ be 
the measure on the $(l-2)$-dimensional sphere 
$\{|\eta^\prime |=(1-|\eta_1|^2)^{\frac{1}{2}}\}$.
If we make a change of variables in the above integral letting 
$\eta^\prime \to (1-|\eta_1|^2)^{\frac{1}{2}}\eta^\prime$, we get
$$
\widehat{B_{z,M}}(y)=
\psi(z)\int_{0}^{+\infty}
r^{-\frac{n}{m}z+l-1}\chi\left(\frac{r}{M}\right)
\int_{-1}^{1}
\frac{(1-|\eta_1|^2)^{\frac{l-2}{2}}}
{|\eta_1|^{\frac{nz}{m}}} \times
$$
$$
\int_{S^{l-2}}F^z(r,\,\eta_1,\,\eta^\prime(1-|\eta_1|^2)^{\frac{l}{2}},\,y) 
d\eta_1d\eta^\prime.
\tag 3.18
$$
We recall that  the distribution
$\dsize \frac{r^{-\frac{n}{m}z+l-1}}{\Gamma(l-\frac{nz}{m})}$ is an entire 
function of $z$ which coincides with the Dirac distribution $\delta_0$ when 
$z =\frac{ml}{n}$, and that the distribution 
$\dsize \frac{|\eta_1|^{\frac{-nz}{m}}}{\Gamma(-\frac{zn}{2m}+\frac{1}{2})}$
is an entire function of $z$ which coincides with the $(k-1)$th derivative of 
the Dirac distribution $\delta_0$ when $l=2k+1$ and $z =\frac{ml}{n}$.
With that in mind  we consider a smooth function $\rho(t)$ which is 
$\equiv 1$ when 
$t\in (0,\,\frac{1}{4})$ and is $\equiv 0$ when 
$t\in (\frac{3}{4},\, 1)$, and we write
\vskip .125 in
$$ 
\widehat{B_{z,M}}(y)=
\psi(z)\left(\int_{-1}^{1}\rho(\eta_1^2) \cdots
+ \int_{-1}^{1}(1-\rho(\eta_1^2))\cdots\right) 
 = \psi(z)I_1(z,M,y)+\psi(z)I_2(z,M,y).
\tag 3.19$$
Since 
$$
(1-\rho(\eta_1^2)) 
(1-\eta_1^2)^{\frac{l-2}{2}}\int_{S^{l-2}}
F^z(r,\,\eta_1,\,\eta^\prime(1-|\eta_1|^2)^{\frac{l}{2}},\,y)d\eta^\prime
d\eta_1 =F_2(r,\,\eta_1)
\tag 3.20
$$
is bounded and continous with respect to $r$ and $\eta_1$, we can write
$$
\psi(z)I_2(z, M,y) = 
\Gamma^{-1}(-\frac{zn}{2m}+\frac{1}{2})\chi\left(\frac{r}{M}\right)
\int_{-1}^{1}F_2(r,\,\eta_1)d\eta_1\Big\vert_{r=0}.
\tag 3.21
$$ 
Thus $|\psi(z)I_2(z,M, y)|$ is bounded by 
a constant with does not depend on $M$ and has at most exponential growth
with respect to Im$(z)$.

We  shall now estimate  $\psi(z)I_1(z, M,y)$.
By our previous observations, the function  
$$
\rho(\eta_1^2) 
(1-\eta_1^2)^{\frac{l-2}{2}}\int_{S^{l-2}}
F^z(r,\,\eta_1,\,\eta^\prime(1-|\eta_1|^2)^{\frac{l}{2}},\,y)d\eta^\prime
d\eta_1= F_1(r,\,\eta_1)
\tag 3.22
$$
can be differentiated $l-1$ times with respect to $\eta_1$ and its derivatives are  continous functions of $r$. Then, if $l= 2k+1$, 
$$
\psi(z)I_1(z, M,y) =
\frac{d^{k-1}}{d\eta_1}\left(\chi\left(\frac{r}{M}\right)F_1(r,\,\eta_1)
\right)\Big\vert_{\eta_1=r=0}.
\tag 3.23
$$
If  $l\ne 2k+1$, and if  $h = [\frac{ \text{ Re}(l)}{2}] $ we use the formula 
$(3)$ in \cite{GS}, pg. $51$, obtaining:
$$
I_1(z,M,y)= \chi\left(\frac r M\right)
\int^{+\infty}_{0}\eta_1^{-\frac{nz}{m}} \Big(
F_1(r,\,\eta_1)+ F_1(r,\, -\eta_1)-
$$
$$
{{2\Big(F_1(r,\,0)+\frac{\eta_1^2}{2!}
\frac{\partial^2}{\partial \eta_1}F_1(r,\, 0) +\cdots + 
\frac{\eta_1^{2h}}{(2h)!}\frac{\partial^{2h}}{\partial \eta_1}
F_1(r,\, 0)\Big)\Big)d\eta_1}_\Big\vert}_{r=0}
\tag 3.24 $$
Thus, $|\psi(z)I_1(z,M,y)|$ is bounded by a  constant which does not depend 
on  $M$. This shows that $\widehat{K_{z,M}}$ is bounded by a uniform constant.  
An easy adaptation of the above argument  shows
that $\widehat{K_{z,M} }$ converges
to $\widehat{K_{z}}$ in distribution sense as $M\to \infty$.
\par
This completes the proof of the theorem.
\enddemo
\vskip.25in 

{\bf Acknowledgements : } The authors wish to thank Eric Sawyer for 
financial support provided through his NSERC grant. The second author wishes to
thank Professor Sawyer for numerous helpful conversations over the last two 
years, which had a profound influence on his work. 

The authors wish to thank Fulvio Ricci for suggesting Theorem 3.1, 
and for related helpful suggestions. 

\newpage
\centerline{\bf References}
\vskip.125in 

\ref \key MC \by M. Christ 
\paper Restriction of the fourier transform to submanifolds of low 
codimension (Ph.D Thesis) U. of Chicago \year 1981 
\endref 

\ref \key DCI \by L. De Carli and A. Iosevich 
\paper A restriction theorem for flat manifolds of codimension two
\yr 1995 \jour Illinois M.J. (to appear) 
\endref 

\ref \key Gr \by A. Greenleaf \pages 519-537
\paper Principal curvature and harmonic analysis 
\yr 1982 \vol30 
\jour Indiana Math J. 
\endref 

\ref \key GS \by I. M. Gel'fand and G. E. Shilov 
\paper Generalized functions vol. 1 
\yr 1964 
\jour  Academic Press.  
\endref

\ref \key H \by H. Hironaka \pages 109-203, 205-326
\paper Resolution of singularities of an algebraic variety over a field of
characteristic $0$ 
\yr 1964 \vol 79
\jour Annals of Math. 
\endref 

\ref \key I1 \by A. Iosevich
\paper Maximal operators associated to families of flat curves in the plane
\yr 1994
\jour Duke Math J. Vol.76
\endref 

\ref \key I2 \by A. Iosevich
\paper Maximal averages over homogeneous hypersurfaces in ${\Bbb R}^3$ 
\yr 1995 \jour Forum Mathematicum (to appear) 
\endref 

\ref \key IS \by A. Iosevich and E. Sawyer 
\paper Oscillatory integrals and maximal averaging operators associated to 
homogeneous hypersurfaces 
\yr 1995 \jour Duke M.J. (to appear)
\endref 

\ref \key L \by W. Littman \pages 766-770
\paper Fourier transforms of surface carried measures and differentiability of
surface averages 
\yr 1963 \vol 69
\jour Bull. Amer. Math. Soc. 
\endref 

\ref \key WR \by W. Rudin
\paper Functional Analysis 
\yr 1973 \jour McGraw-Hill Inc. 
\endref 

\ref \key P \by E. Prestini 
\paper Restriction theorems for the fourier transform to some manifolds in 
${\Bbb R}^n$  
\yr 1979 \jour Proceedings of Symposia in Pure Mathematics 
\vol XXXV Part 1 
\endref 

\ref \key TS \by T. Shintani \pages 25-65
\paper On zeta function associated with the vector space of quadratic forms
\yr 1975 \jour J. of F.Sc. U. of Tokyo
\endref

\ref \key So  \by C. D. Sogge         
\paper Fourier Integrals in Classical Analysis 
\yr1991 \jour Oxford University Press 
\endref 

\ref \key St  \by E. M. Stein 
\paper Harmonic Analysis 
\yr 1993 \vol43
\jour Princeton University Press 
\endref 

\enddocument